\documentclass[preprint,review,12pt]{elsarticle}
\usepackage{amsmath,amssymb,amsthm}
\usepackage{booktabs,multirow,graphicx,subcaption}
\usepackage{algorithm,algpseudocode}
\usepackage{xcolor,placeins}
\usepackage[hidelinks]{hyperref}
\biboptions{sort&compress}

\newcommand{\R}{\mathbb{R}}
\newcommand{\bd}{\partial}
\newcommand{\D}{\Delta}
\newcommand{\LL}{\mathcal{L}}
\newcommand{\xx}{\mathbf{x}}
\newcommand{\yy}{\mathbf{y}}
\newcommand{\vv}{\mathbf{v}}
\newcommand{\pp}{\mathbf{p}}
\newcommand{\ww}{\mathbf{w}}
\newcommand{\ba}{\mathbf{a}}
\newcommand{\rev}[1]{{\color{black}#1}}
\hypersetup{pdfauthor={Huifang Zhou},pdftitle={Boundary-preserving randomized neural network methods for Poisson equations on general two- and three-dimensional domains}}
\begin{document}
\begin{frontmatter}
\title{Boundary-preserving randomized neural network methods for Poisson equations on general two- and three-dimensional domains}
\author[aff1]{Huifang Zhou\corref{cor1}}
\ead{zhouhuifang@jlu.edu.cn}
\cortext[cor1]{Corresponding author.}
\affiliation[aff1]{organization={School of Mathematics, Jilin University},city={Changchun},postcode={130012},country={P.R. China}}
\begin{abstract}
In this paper, we propose a boundary-preserving randomized neural network method for Poisson equations with Dirichlet boundary conditions in two and three dimensions. A linear lifting imposes the boundary conditions exactly, without balancing penalty weights between interior and boundary residuals. We construct the liftings by generalized barycentric coordinates and transfinite interpolation on convex and nonconvex polygons and polyhedra, and by radial parameterization on star-shaped curved domains. Numerical experiments on different domains show that the proposed method substantially reduces errors compared with the standard randomized neural network method.
\end{abstract}
\end{frontmatter}

\clearpage
\section{Introduction}
\label{sec:intro}

\rev{In recent years, deep learning methods have been widely used to solve partial differential equations (PDEs). Compared with grid-based methods, these methods can avoid mesh generation and provide flexible approximations for high-dimensional PDEs \cite{sirignano2018}.} \rev{Physics-informed neural networks (PINNs) incorporate the equation residuals together with boundary and initial conditions into the loss function and can handle both forward and inverse problems \cite{raissi2019,wang2021,karniadakis2021}. The Deep Galerkin method minimizes the strong-form PDE residual at randomly sampled collocation points \cite{sirignano2018,alaradi2022,jiang2026}. For second-order elliptic problems with a variational structure, the Deep Ritz method minimizes the corresponding energy functional, which involves only first-order derivatives \cite{eyu2018,lu2021,muller2022}. These methods turn the solution of a PDE into an optimization problem over the network parameters, which are typically optimized by stochastic gradient methods.} When the hidden-layer parameters are trained, the optimization is generally nonlinear, even for a linear PDE.

A randomized neural network (RNN) fixes the hidden-layer parameters after random sampling and determines only the output coefficients \rev{\cite{igelnik1995,huang2006}}. The hidden-layer parameters require no iterative training, which saves the associated optimization cost. For PDEs, \rev{the RNN method} can be combined with domain decomposition \cite{dong2021}. \rev{The sampling range of the random hidden-layer parameters can also be optimized \cite{dong2022}.} Multiscale random features and residual rescaling provide another collocation approach \cite{chen2022}. \rev{Related methods have been developed for elliptic equations with sharp gradients, time-dependent equations, and interface problems \cite{calabro2021,chen2023,chi2024}. However, approximate enforcement of boundary conditions can introduce additional errors.}

\rev{Dirichlet boundary conditions in neural network methods can be imposed through the training loss or built into the trial function.} \rev{Standard PINNs impose soft constraints by including the boundary and interior residuals in a weighted loss \cite{raissi2019}.} \rev{An adaptive weighting strategy adjusts the weights of the PDE and boundary losses using gradient information to reduce their imbalance during training \cite{wang2021}.} \rev{The Deep Ritz method analyzed in \cite{muller2022} adds a boundary penalty term to the energy functional.} \rev{The Deep Nitsche method imposes the boundary conditions weakly by adding consistency and boundary-penalty terms to the variational functional \cite{liao2021}.} \rev{These approaches use penalty or stabilization parameters \cite{raissi2019,wang2021,muller2022,liao2021}. In \cite{berg2018,sheng2021}, a network fitted to the boundary data is combined with an interior network through a distance function or a length factor. Hard constraints construct trial functions that satisfy the boundary conditions exactly \cite{lagaris1998}. Length factors allow exact enforcement of mixed Dirichlet and Neumann conditions on irregular boundaries \cite{mcfall2009}. Methods based on the theory of functional connections incorporate boundary constraints analytically into trial functions involving neural networks \cite{leake2020,schiassi2021}. Approximate distance functions and transfinite interpolation provide such constructions on complex domains \cite{sukumar2022}, thereby avoiding boundary penalty weights.} A comparison of four \rev{boundary} enforcement strategies shows that approximate distance functions can improve accuracy~\cite{berrone2023}. \rev{For elliptic equations on rectangular domains, boundary interpolation is combined with the RNN method to impose Dirichlet or clamped boundary conditions exactly \cite{zhou2026}.}

To impose Dirichlet conditions by a trial function on a general domain, we need an interpolation operator that reproduces the boundary trace and has computable derivatives. On polygons and polyhedra, Wachspress coordinates \cite{wachspress1975} and mean value coordinates \cite{floater2003,hormann2006,ju2005} reproduce affine functions from vertex values. Transfinite interpolation incorporates \rev{values along edges and across faces} to reproduce the full trace \cite{randrianarivony2011}. For star-shaped curved domains with a positive radial graph, polar parameterization connects an interior center to the boundary along rays \cite{juttler2019}. These constructions give two ways to define the lifting used below.

\rev{In this paper, we combine these geometric constructions to give a unified framework for boundary-preserving (BP) RNN methods on general two- and three-dimensional domains. We use a linear lifting operator to construct neural network approximations that satisfy Dirichlet boundary conditions exactly. For polygonal and polyhedral domains, we construct barycentric--transfinite lifting operators using generalized barycentric coordinates and transfinite interpolation. We use Wachspress coordinates for convex polytopes and mean value coordinates for nonconvex polytopes. For star-shaped curved domains, we construct radial lifting operators using radial parameterization. In the reported numerical tests, the BP method reduces the errors relative to the standard RNN method by an average factor of approximately $190$, with a maximum of approximately $1000$.}

The paper is organized as follows. Section~\ref{sec:framework} presents the formulation of the BP method. Section~\ref{sec:poly} constructs barycentric--transfinite liftings on polygons and polyhedra. Section~\ref{sec:radial} gives radial liftings on star-shaped curved domains. Numerical experiments are reported in Section~\ref{sec:numerics}, and conclusions are given in Section~\ref{sec:conclusion}.

\section{Problem setting and the boundary-preserving framework}
\label{sec:framework}

In this section, we state the Poisson problem and introduce the standard RNN method with boundary residuals. We then construct the trial function of the BP method and derive a linear least-squares system using only interior equations.

\subsection{\rev{Poisson problem and the standard RNN method}}

Let $\Omega\subset\R^d$, $d=2,3$, be a bounded domain with boundary
$\bd\Omega$. In this paper, we consider the Poisson problem
\begin{equation}
\label{eq:poisson}
-\D u(\xx)=f(\xx)\quad\text{in }\Omega,
\qquad
u(\xx)=g(\xx)\quad\text{on }\bd\Omega,
\end{equation}
where $f:\Omega\to\R$ and $g:\bd\Omega\to\R$ are the given source and boundary data.

We approximate the solution by the RNN
\begin{equation*}
\hat u(\xx;\ww)=\sum_{j=1}^{M} w_j\,\phi_j(\xx),
\qquad
\phi_j(\xx)=\tanh(\ba_j\cdot\xx+b_j),
\end{equation*}
where the vectors $\ba_j\in\R^d$ and scalars $b_j\in\R$ are sampled from a fixed
distribution and remain fixed. Only the output weights $\ww=(w_1,\dots,w_M)^\top$ are
unknown. This is a single-hidden-layer random feature representation
\cite{rahimi2007}.

The standard RNN method imposes the PDE and boundary conditions through residual
equations. We collocate the
PDE at interior points and the boundary condition at boundary points. Then the
output weights are determined by solving the following linear system in the
least-squares sense:
\begin{equation*}
-\sum_{j=1}^{M} w_j\,\D\phi_j(\xx_i)=f(\xx_i),
\qquad i=1,\dots,N_{\mathrm{int}},
\end{equation*}
\begin{equation}
\label{eq:std_bnd}
\sum_{j=1}^{M} w_j\,\phi_j(\yy_k)=g(\yy_k),
\qquad k=1,\dots,N_{\mathrm{bnd}},
\end{equation}
where $N_{\mathrm{int}}$ and $N_{\mathrm{bnd}}$ denote the numbers of interior and boundary collocation points, respectively. The points satisfy $\xx_i\in\Omega$ and $\yy_k\in\bd\Omega$. The boundary condition is imposed
through \eqref{eq:std_bnd}, so the relation
$\hat u=g$ on $\bd\Omega$ is not satisfied exactly. We call this the soft-boundary
formulation.

\subsection{Boundary-preserving formulation}

Let $\LL$ be a linear lifting that extends a function $h$ on $\bd\Omega$
to $\Omega$ and satisfies the boundary interpolation property
\begin{equation}
\label{eq:interp}
(\LL h)|_{\bd\Omega}=h|_{\bd\Omega}.
\end{equation}
We construct the boundary-preserving trial function as
\begin{equation}
\label{eq:bp}
\hat u(\xx;\ww)=\LL g(\xx)+\sum_{j=1}^{M} w_j\bigl(\phi_j(\xx)-\LL\phi_j(\xx)\bigr).
\end{equation}

Here $\LL g$ gives the prescribed boundary values, and each corrected feature $\phi_j-\LL\phi_j$ vanishes on the boundary. The correction therefore changes only the interior approximation.

By \eqref{eq:interp}, $\LL g=g$ and $\phi_j-\LL\phi_j=0$ on $\bd\Omega$, so
\begin{equation*}
\hat u|_{\bd\Omega}=g
\end{equation*}
holds for every weight vector $\ww$. Thus, the boundary condition is satisfied by
construction, and no boundary penalty needs to be balanced against the interior
residual. We determine the weights by collocating only the interior PDE. Since
\begin{equation*}
\D\hat u=\D\LL g+\sum_{j=1}^{M}w_j\bigl(\D\phi_j-\D\LL\phi_j\bigr),
\end{equation*}
substituting \eqref{eq:bp} into \eqref{eq:poisson} yields the linear system
\begin{equation*}
-\sum_{j=1}^{M} w_j\Bigl[\D\phi_j(\xx_i)-\D\LL\phi_j(\xx_i)\Bigr]
= f(\xx_i)+\D\LL g(\xx_i),
\qquad i=1,\dots,N_{\mathrm{int}},
\end{equation*}
which is solved in the least-squares sense.

We also consider affine reproduction: $\LL p=p$ for every affine function $p$. This property is useful but is not required for boundary preservation. We next construct $\LL$ for two classes of domains. Section~\ref{sec:poly} combines generalized barycentric coordinates with transfinite interpolation on polytopes. Section~\ref{sec:radial} gives a radial construction for star-shaped curved domains.

The solution procedure is summarized in Algorithm~\ref{alg:bp}.

\begin{algorithm}[htbp]

\captionsetup{labelfont={color=black},textfont={color=black}}
\caption{Boundary-preserving randomized neural network method}
\label{alg:bp}
\begin{algorithmic}[1]
\Require Domain $\Omega$, data $f,g$, feature count $M$, interior points $\{\xx_i\}_{i=1}^{N_{\mathrm{int}}}$
\Ensure Approximation $\hat u(\xx;\ww^\star)$
\State Construct the lifting $\LL$ under the geometry and regularity conditions in Sections~\ref{sec:poly} and \ref{sec:radial}.
\State Sample and fix $\{\phi_j\}_{j=1}^M$; assemble $A_{ij}=-\D(\phi_j-\LL\phi_j)(\xx_i)$ and $b_i=f(\xx_i)+\D\LL g(\xx_i)$.
\State Solve the least-squares problem: $\displaystyle\ww^\star\in\operatorname*{arg\,min}\limits_{\ww\in\R^M}\|A\ww-\mathbf b\|_2^2$.
\State \Return $\hat u(\xx;\ww^\star)$.
\end{algorithmic}
\end{algorithm}
\FloatBarrier

\section{Barycentric liftings on polygonal and polyhedral domains}
\label{sec:poly}

In this section, we construct the lifting \rev{operator} for polygons and polyhedra. We use generalized barycentric coordinates to express a point as a weighted combination of the vertices. Wachspress coordinates are used for convex polytopes and mean value coordinates for nonconvex polytopes. \rev{Finally, we combine these coordinates with transfinite interpolation to obtain linear lifting operators.}

\subsection{Generalized barycentric coordinates}

Let $\mathcal P$ be a polygon in $\R^2$ or a polyhedron in $\R^3$ with vertices
$\vv_1,\dots,\vv_n$. We assume that $\mathcal P$ is bounded. Its boundary is a
simple closed polygon in $\R^2$, which may be nonconvex but has no
self-intersections or holes, or a closed two-manifold polyhedral surface in
$\R^3$, which may be nonconvex and may have holes or cavities. Neither case
requires a star-shaped domain. A family of generalized barycentric coordinates
$\lambda_1,\dots,\lambda_n$ on $\mathcal P$ satisfies
\begin{equation}
\label{eq:coords}
\sum_{i=1}^{n}\lambda_i(\xx)=1,
\qquad
\sum_{i=1}^{n}\lambda_i(\xx)\vv_i=\xx,
\end{equation}
and $\lambda_i$ vanishes on each codimension-one face of $\mathcal P$ that does not
contain $\vv_i$. Here a codimension-one face is an edge in two dimensions and a
face in three dimensions. The partition of unity and linear reproduction in
\eqref{eq:coords} imply that the coordinates reproduce affine functions. For every affine function $p$, the sum $\sum_i\lambda_i(\xx)p(\vv_i)$ equals $p(\xx)$. The constant part of $p$ is reproduced by $\sum_i\lambda_i(\xx)=1$, and the linear part by $\sum_i\lambda_i(\xx)\vv_i=\xx$.

For a convex polygon, the Wachspress coordinates \cite{wachspress1975} are
defined by the weights
\begin{equation*}
\beta_i(\xx)=\frac{A(\vv_{i-1},\vv_i,\vv_{i+1})}
{A(\vv_{i-1},\vv_i,\xx)\,A(\vv_i,\vv_{i+1},\xx)},
\qquad
\lambda_i(\xx)=\frac{\beta_i(\xx)}{\sum_{j=1}^{n}\beta_j(\xx)},
\end{equation*}
where $A(\mathbf a,\mathbf b,\mathbf c)=\tfrac12\det(\mathbf b-\mathbf a,\mathbf c-\mathbf a)$ denotes the signed area of the triangle with vertices $\mathbf a,\mathbf b,\mathbf c$; it is positive for counterclockwise orientation. We number the polygon vertices counterclockwise. The coordinates are nonnegative on a convex polygon.

For a convex polyhedron, let $F_1,\dots,F_k$ be the facets incident to $\vv_i$, ordered cyclically. Their outward normals form the rays of the normal cone at $\vv_i$. A fan decomposition into triples $(F_j,F_{j+1},F_k)$ gives the weight
\begin{equation*}
\beta_i(\xx)=\sum_{j=1}^{k-2}
\frac{\bigl|\det(\mathbf{n}_{F_j},\mathbf{n}_{F_{j+1}},\mathbf{n}_{F_k})\bigr|}
{d_{F_j}(\xx)\,d_{F_{j+1}}(\xx)\,d_{F_k}(\xx)},
\end{equation*}
where $\mathbf{n}_F$ denotes the outward unit normal of facet $F$, and $d_F(\xx)$ denotes the
perpendicular distance from $\xx$ to the supporting plane of $F$.

For a nonconvex polygon, we use mean value coordinates \cite{floater2003,hormann2006}, since the nonnegativity and nonsingularity of Wachspress coordinates on convex polygons do not extend to this case. The weights are
\begin{equation*}
\beta_i(\xx)=\frac{\tan\bigl(\alpha_{i-1}(\xx)/2\bigr)
+\tan\bigl(\alpha_i(\xx)/2\bigr)}{\|\vv_i-\xx\|},
\qquad
\lambda_i(\xx)=\frac{\beta_i(\xx)}{\sum_{j=1}^{n}\beta_j(\xx)},
\end{equation*}
where $\alpha_i(\xx)\in(-\pi,\pi)$ denotes the oriented angle from $\vv_i-\xx$ to $\vv_{i+1}-\xx$.

For a nonconvex polyhedron, we use the three-dimensional mean value coordinates
of Ju--Schaefer--Warren \cite{ju2005}. The boundary is represented by a closed
triangular mesh $\mathcal T$ with consistent outward orientation. Let $\xx$ be an
interior point. For each triangle $T=(\vv_i,\vv_j,\vv_k)\in\mathcal T$, we first compute its contributions to the three vertex weights. Define the unit vectors from $\xx$ toward these vertices by
\begin{equation*}
\mathbf e_{i}=\frac{\vv_{i}-\xx}{\|\vv_{i}-\xx\|},\quad
\mathbf e_j=\frac{\vv_j-\xx}{\|\vv_j-\xx\|},\quad
\mathbf e_k=\frac{\vv_k-\xx}{\|\vv_k-\xx\|}.
\end{equation*}
The corresponding spherical triangle has side lengths $\theta_i=\arccos(\mathbf e_j\cdot\mathbf e_k)$, $\theta_j=\arccos(\mathbf e_k\cdot\mathbf e_i)$, and $\theta_k=\arccos(\mathbf e_i\cdot\mathbf e_j)$.
With $\sigma=(\theta_{i}+\theta_{j}+\theta_{k})/2$, the quantities associated with $\vv_i$ are
\begin{equation*}
c_{i}=\frac{2\sin\sigma\,\sin(\sigma-\theta_{i})}{\sin\theta_{j}\,\sin\theta_{k}}-1,
\qquad
s_{i}=\operatorname{sgn}\bigl(\det[\mathbf e_{i},\mathbf e_{j},\mathbf e_{k}]\bigr)\sqrt{1-c_{i}^2}.
\end{equation*}
The remaining quantities follow by cyclic permutation of $(i,j,k)$. The contribution of the triangle $T$ to the weight of vertex $\vv_{i}$ is
\begin{equation*}
\mu_{i}^T=\frac{\theta_{i}-c_{j}\theta_{k}-c_{k}\theta_{j}}
{\|\vv_{i}-\xx\|\,\sin\theta_{j}\,s_{k}}.
\end{equation*}
The same index replacements give $\mu_j^T$ and $\mu_k^T$. We sum the contributions at each vertex $\vv_i$ over all triangles containing $\vv_i$ and normalize the weights:
\begin{equation*}
\beta_i(\xx)=\sum_{T\in\mathcal T:\,\vv_i\in T}\mu_i^T,
\qquad
\lambda_i(\xx)=\frac{\beta_i(\xx)}{\sum_{j=1}^{n}\beta_j(\xx)}.
\end{equation*}

On a convex polytope, the nonnegative coordinates satisfy $0\le\lambda_i\le1$ and sum to one. Thus, $\sum_i\lambda_i h(\vv_i)$ is a convex combination of the vertex values and lies between their minimum and maximum. The edge and face evaluation points defined below also remain on the finite boundary pieces, so no extension is needed. On a nonconvex polytope, mean value coordinates may be negative or exceed one. The evaluation points may then leave those boundary pieces, which is why the lifting requires the boundary-data extension in Section~\ref{sec:extension}.

\subsection{Two-dimensional transfinite lifting}

We number the polygon vertices cyclically. With $\vv_i$ as the base vertex, define one evaluation point on each of its two incident edge lines:
\begin{equation*}
\begin{aligned}
\pp_i^{+}(\xx)&=\vv_i+\lambda_{i+1}(\xx)(\vv_{i+1}-\vv_i),\\
\pp_i^{-}(\xx)&=\vv_i+\lambda_{i-1}(\xx)(\vv_{i-1}-\vv_i).
\end{aligned}
\end{equation*}
The superscripts $+$ and $-$ only distinguish the two edges at the same base vertex; they do not indicate signs or arithmetic operations. On the edge $[\vv_i,\vv_{i+1}]$, $\pp_i^+(\xx)=\pp_{i+1}^-(\xx)=\xx$. Set
\begin{equation*}
S_i[h](\xx)=h(\pp_i^{+}(\xx))+h(\pp_i^{-}(\xx))-h(\vv_i),
\end{equation*}
with the transfinite lifting given by
\begin{equation}
\label{eq:2dtrans}
\LL h(\xx)=\sum_{i=1}^{n}\lambda_i(\xx)\,S_i[h](\xx).
\end{equation}
On each edge, \eqref{eq:2dtrans} reproduces $h$ by the boundary interpolation property
of $\lambda_i$. Thus, \eqref{eq:interp} holds. Applying the chain rule twice to
\eqref{eq:2dtrans} gives the Laplacian

\begin{align*}
\D\LL h(\xx)&=\sum_{i=1}^{n}\Bigl[\D\lambda_i(\xx)\,S_i[h](\xx)
+2\nabla\lambda_i(\xx)\cdot\nabla S_i[h](\xx)
+\lambda_i(\xx)\,\D S_i[h](\xx)\Bigr],
\end{align*}
where all indices are cyclic and
\begin{align}
\nabla S_i[h](\xx)={}&
\bigl[(\mathbf d_i\cdot\nabla)h(\pp_i^{+}(\xx))\bigr]\nabla\lambda_{i+1}(\xx)\nonumber\\
&-\bigl[(\mathbf d_{i-1}\cdot\nabla)h(\pp_i^{-}(\xx))\bigr]\nabla\lambda_{i-1}(\xx),\label{eq:2dgrad}\\[2pt]
\D S_i[h](\xx)={}&
\bigl[(\mathbf d_i\cdot\nabla)h(\pp_i^{+}(\xx))\bigr]\D\lambda_{i+1}(\xx)\nonumber\\
&+\Bigl[(\mathbf d_i\cdot\nabla)^2h(\pp_i^{+}(\xx))\Bigr]
\bigl|\nabla\lambda_{i+1}(\xx)\bigr|^2\nonumber\\
&-\bigl[(\mathbf d_{i-1}\cdot\nabla)h(\pp_i^{-}(\xx))\bigr]\D\lambda_{i-1}(\xx)\nonumber\\
&+\Bigl[(\mathbf d_{i-1}\cdot\nabla)^2h(\pp_i^{-}(\xx))\Bigr]
\bigl|\nabla\lambda_{i-1}(\xx)\bigr|^2.\label{eq:2dlapSi}
\end{align}
Here $\mathbf d_i=\vv_{i+1}-\vv_i$ and $\mathbf d_{i-1}=\vv_i-\vv_{i-1}$. For $h=g$, the directional derivatives are the first and second tangential derivatives of the boundary trace along the corresponding edges. For $h=\phi_j$, we differentiate the analytic $\tanh$ function.

\subsection{Three-dimensional transfinite lifting}

For a polyhedron, we use the transfinite blending formula of Randrianarivony
\cite{randrianarivony2011}. For each vertex $\vv_i$, let $\mathcal F_i$ contain all
faces incident to $\vv_i$. Here ``face'' includes the vertex itself
(dimension $0$), an edge (dimension $1$), and a two-dimensional face
(dimension $2$). For $\Pi\in\mathcal F_i$, define
\begin{equation*}
\pp_{\Pi,i}(\xx)=\vv_i+\sum_{\substack{\vv_j\in\Pi\\j\ne i}}\lambda_j(\xx)\,
\bigl(\vv_j-\vv_i\bigr).
\end{equation*}
The sum runs over the vertices of $\Pi$ other than $\vv_i$. For the vertex face $\Pi=\{\vv_i\}$, the sum is empty and $\pp_{\Pi,i}=\vv_i$. Set
\begin{equation*}
S_i[h](\xx)=\sum_{\Pi\in\mathcal F_i}(-1)^{2-\dim\Pi}h(\pp_{\Pi,i}(\xx)),
\end{equation*}
where $\dim\Pi\in\{0,1,2\}$ denotes the topological dimension of
$\Pi$ ($0$ for a vertex, $1$ for an edge, and $2$ for a face). Thus, $S_i[h]$
contains the vertex term $h(\vv_i)$ with sign $+1$, the edge terms with sign
$-1$, and the face terms with sign $+1$. We define the lifting by
\begin{equation}
\label{eq:3dtrans}
\LL h(\xx)=\sum_{i=1}^{n}\lambda_i(\xx)\,S_i[h](\xx).
\end{equation}
We use each geometric face without holes directly as a two-dimensional face in the formula. For a face with holes, we first subdivide it conformingly into patches without holes and use each patch as a two-dimensional face in the formula. To reproduce the boundary trace, the resulting faces must be compatible. Each face is a disk, every edge belongs to two faces, and the faces around each vertex form a cycle. Two faces may intersect only along a common edge or at a common vertex. The coordinates of vertices outside a boundary face must vanish on that face. Under these conditions, the alternating sum reduces to $h(\xx)$ on the boundary.

Applying the product rule and the chain rule to \eqref{eq:3dtrans} gives
\begin{equation*}
\D\LL h(\xx)=\sum_{i=1}^{n}\Bigl[\D\lambda_i(\xx)\,S_i[h](\xx)
+2\nabla\lambda_i(\xx)\cdot\nabla S_i[h](\xx)
+\lambda_i(\xx)\,\D S_i[h](\xx)\Bigr].
\end{equation*}
For each $\Pi\in\mathcal F_i$, choose an orthonormal tangent basis $\{\mathbf t_\alpha\}_{\alpha=1}^{\dim\Pi}$. There is one vector for an edge, two for a face, and none for a vertex. Then

\begin{align}
\nabla S_i[h](\xx)&=\sum_{\Pi\in\mathcal F_i}(-1)^{2-\dim\Pi}\sum_{\alpha}
\bigl[(\mathbf t_{\alpha}\cdot\nabla)h(\pp_{\Pi,i}(\xx))\bigr]\,J_{\Pi,i}(\xx)^\top \mathbf t_{\alpha},\label{eq:3dgrad}\\
\D S_i[h](\xx)&=\sum_{\Pi\in\mathcal F_i}(-1)^{2-\dim\Pi}
\Bigl[\sum_{\alpha,\beta}\bigl[(\mathbf t_{\alpha}\cdot\nabla)(\mathbf t_{\beta}\cdot\nabla)h(\pp_{\Pi,i}(\xx))\bigr]
\nonumber\\
&\qquad\qquad\times\bigl(J_{\Pi,i}(\xx)^\top \mathbf t_{\alpha}\bigr)\cdot\bigl(J_{\Pi,i}(\xx)^\top \mathbf t_{\beta}\bigr)\nonumber\\
&\qquad+\sum_{\alpha}\bigl[(\mathbf t_{\alpha}\cdot\nabla)h(\pp_{\Pi,i}(\xx))\bigr]\bigl(\mathbf t_{\alpha}\cdot\D\pp_{\Pi,i}(\xx)\bigr)\Bigr],\label{eq:3dlapSi}
\end{align}
where $J_{\Pi,i}(\xx)=D_{\xx}\pp_{\Pi,i}(\xx)=\sum_{\substack{\vv_j\in\Pi\\j\ne i}}(\vv_j-\vv_i)\nabla\lambda_j(\xx)^\top$ is the Jacobian of the evaluation map.
We compute the Wachspress coordinate derivatives analytically. For mean value coordinates, the angles and arc lengths make analytic differentiation inconvenient. We therefore use second-order centered finite differences with step $1\times10^{-5}$ in both two and three dimensions.

\subsection{Boundary-data extension for nonconvex domains}
\label{sec:extension}

On a nonconvex polytope, signed mean value coordinates can place the evaluation points used to compute $\LL h$ outside $\partial\Omega$, while they remain on the corresponding edge lines or face planes. We therefore extend the boundary data to these lines and planes using quartic Taylor polynomials.

In two dimensions, assume that $g\in C^4$ on each edge. Write $g_e(t)=g(\vv_i+t(\vv_{i+1}-\vv_i))$, $t\in[0,1]$, on $[\vv_i,\vv_{i+1}]$, and extend it to the whole edge line by
\begin{equation*}
\tilde g_e(t)=
\begin{cases}
\sum_{k=0}^{4}\frac{g_e^{(k)}(0)}{k!}\,t^k,&t<0,\\[2pt]
g_e(t),&0\le t\le1,\\[2pt]
\sum_{k=0}^{4}\frac{g_e^{(k)}(1)}{k!}\,(t-1)^k,&t>1
\end{cases},
\end{equation*}
where $g_e^{(k)}$ is the $k$-th derivative along the edge and is determined by
the boundary trace $g|_{\partial\Omega}$. The two terms in $S_i[g]$ use
$\lambda_{i+1}(\xx)$ and $1-\lambda_{i-1}(\xx)$ as parameters. Since
$\pp_i^+(\xx)=\vv_i+\lambda_{i+1}(\xx)\mathbf d_i$, we replace $g$ by
$\tilde g_e$ along this edge line in \eqref{eq:2dgrad} and \eqref{eq:2dlapSi}, giving
\begin{equation*}
(\mathbf d_i\cdot\nabla)g(\pp_i^{+}(\xx))\;\mapsto\;\tilde g_e'(\lambda_{i+1}(\xx)),
\qquad
(\mathbf d_i\cdot\nabla)^2g(\pp_i^{+}(\xx))\;\mapsto\;\tilde g_e''(\lambda_{i+1}(\xx)).
\end{equation*}
For $\pp_i^-(\xx)$, we use the parameter $1-\lambda_{i-1}(\xx)$ on the edge
$[\vv_{i-1},\vv_i]$ in the same way.

In three dimensions, we assume convex faces with facewise traces $g\in C^6$. Vertex terms need no extension, and edge traces use the same quartic endpoint extension as in two dimensions. For a face term, we use $g$ on the finite face and extend it within the face plane. Outside the face, let $\mathbf q$ be the closest point on its boundary. Choose orthonormal in-plane directions $\mathbf n$, $\mathbf t$, normal and tangent to the edge containing $\mathbf q$. At a vertex, choose an incident edge according to a fixed convention. Write
\begin{equation*}
\pp_{\Pi,i}(\xx)=\mathbf q+s\mathbf n+\tau\mathbf t,
\qquad
s=(\pp_{\Pi,i}(\xx)-\mathbf q)\cdot\mathbf n,\quad
\tau=(\pp_{\Pi,i}(\xx)-\mathbf q)\cdot\mathbf t,
\end{equation*}
and use the fourth-order Taylor polynomial in the face plane,
\begin{equation*}
\tilde g(\pp_{\Pi,i}(\xx))
=\sum_{\substack{k,l\ge0\\k+l\le4}}\frac{(\mathbf n\cdot\nabla)^{k}(\mathbf t\cdot\nabla)^{l} g(\mathbf q)}
{k!\,l!}\;s^{k}\tau^{l};
\end{equation*}
we then replace the in-plane tangential derivatives of the face term in
\eqref{eq:3dgrad} and \eqref{eq:3dlapSi} by the derivatives of $\tilde g$:
\begin{equation*}
\begin{aligned}
(\mathbf t_{\alpha}\cdot\nabla)g(\pp_{\Pi,i}(\xx))&\;\mapsto\;(\mathbf t_{\alpha}\cdot\nabla)\tilde g(\pp_{\Pi,i}(\xx)),\\[2pt]
(\mathbf t_{\alpha}\cdot\nabla)(\mathbf t_{\beta}\cdot\nabla)g(\pp_{\Pi,i}(\xx))&\;\mapsto\;(\mathbf t_{\alpha}\cdot\nabla)(\mathbf t_{\beta}\cdot\nabla)\tilde g(\pp_{\Pi,i}(\xx)).
\end{aligned}
\end{equation*}
Here $k,l$ are nonnegative derivative orders, whereas $\mathbf t_\alpha,\mathbf t_\beta$ are the face tangent vectors in \eqref{eq:3dgrad} and \eqref{eq:3dlapSi}, with $\alpha,\beta\in\{1,2\}$. We differentiate the extension within each smooth piece.

For a convex polytope, the Wachspress coordinates belong to $[0,1]$. The evaluation points remain on the finite faces, so no extension is needed. The radial lifting in Section~\ref{sec:radial} also uses only true boundary points and does not require this extension.

\section{Radial liftings on star-shaped curved domains}
\label{sec:radial}

\rev{In this section, we construct radial lifting operators for star-shaped curved domains in two and three dimensions.}

\subsection{Two-dimensional radial lifting}

A \rev{domain} $\Omega\subset\R^2$ is \emph{star-shaped} if it contains an interior point $\mathbf c\in\Omega$ such that, for every $\xx\in\Omega$, the segment $(1-t)\mathbf c+t\xx$, $t\in[0,1]$, remains in $\Omega$. We further assume that the boundary has a positive continuous radius $r(\theta)$ and that each ray meets it once. Then
\begin{equation*}
\bd\Omega=\{\boldsymbol\gamma(\theta)=\mathbf c+r(\theta)\mathbf e_r(\theta)
:\theta\in[0,2\pi)\},
\qquad
\mathbf e_r(\theta)=(\cos\theta,\sin\theta).
\end{equation*}
For $\xx=\mathbf c+\rho\,\mathbf e_r(\theta)$ with $0\le\rho<r(\theta)$, let
$\tilde h(\theta)=h(\boldsymbol\gamma(\theta))$. Choose a fixed reference angle
$\theta_0$ and set $\tilde h_0=h(\boldsymbol\gamma(\theta_0))$. We define the radial lifting by
\begin{equation}
\label{eq:2drad}
\LL h(\xx)=\tilde h_0+\rho^2 q(\theta),
\qquad
q(\theta)=\frac{\tilde h(\theta)-\tilde h_0}{r(\theta)^2}.
\end{equation}
At $\rho=r(\theta)$, \eqref{eq:2drad} gives
$\LL h=\tilde h(\theta)=h(\boldsymbol\gamma(\theta))$. Hence \eqref{eq:interp} holds, and
$\LL$ leaves constants unchanged. Unlike the barycentric liftings in Section~\ref{sec:poly}, this radial lifting generally does not reproduce affine functions. It is quadratic in the radial variable $\rho$, whereas an affine function generally varies linearly along a ray.
Taking the Laplacian of
\eqref{eq:2drad} in polar coordinates gives
$\D\LL h=4q(\theta)+q''(\theta)$.
For an explicit smooth radius, we compute $q''$ in closed form by differentiating
$q$ with respect to $\theta$. When a boundary is given parametrically and its radial function is not available explicitly, the radius is tabulated numerically and its geometric
derivatives are computed numerically.

\subsection{Three-dimensional radial lifting}

Let $\Omega\subset\R^3$ be star-shaped about $\mathbf c$, with radial boundary
\begin{equation*}
\bd\Omega=\{\boldsymbol\gamma(\theta,\varphi)
=\mathbf c+r(\theta,\varphi)\mathbf e_r(\theta,\varphi)\},
\qquad
\mathbf e_r=(\sin\theta\cos\varphi,\sin\theta\sin\varphi,\cos\theta).
\end{equation*}
For $\xx=\mathbf c+\rho\mathbf e_r(\theta,\varphi)$ with
$0\le\rho<r(\theta,\varphi)$, define
$\tilde h(\theta,\varphi)=h(\boldsymbol\gamma(\theta,\varphi))$. For fixed reference
angles $(\theta_0,\varphi_0)$, let
$\tilde h_0=h(\boldsymbol\gamma(\theta_0,\varphi_0))$. We construct the lifting as
\begin{equation}
\label{eq:3drad}
\LL h(\xx)=\tilde h_0+\rho^2 q(\theta,\varphi),
\qquad
q(\theta,\varphi)=\frac{\tilde h(\theta,\varphi)-\tilde h_0}{r(\theta,\varphi)^2}.
\end{equation}
Taking the Laplacian of \eqref{eq:3drad} in spherical coordinates gives
\begin{equation*}
\D\LL h=6q+\D_{S^2}q,
\qquad
\D_{S^2}q=\frac{1}{\sin\theta}\frac{\partial}{\partial\theta}
\Bigl(\sin\theta\,\frac{\partial q}{\partial\theta}\Bigr)
+\frac{1}{\sin^2\theta}\frac{\partial^2 q}{\partial\varphi^2},
\end{equation*}
where $\D_{S^2}$ denotes the Laplace--Beltrami operator on the unit sphere.

\section{Numerical experiments}
\label{sec:numerics}

In this section, we compare the accuracy of the boundary-preserving method with that of the standard RNN method on polygonal, polyhedral, and star-shaped curved domains. We first give the experimental setup, and then present the results in two and three dimensions.

\begin{samepage}
\subsection{Experimental setup}
\label{sec:setup}

We solve the Poisson problem using the exact solution
\begin{equation*}
u(\xx)=\sin(\pi x)\sin(\pi y)
\end{equation*}
in two dimensions and
\begin{equation*}
u(\xx)=\sin(\pi x)\sin(\pi y)\sin(\pi z)
\end{equation*}
in three dimensions. The corresponding source terms are $f=2\pi^2u$ and $f=3\pi^2u$,
respectively.
\end{samepage}

The RNN has a single hidden layer with $M$ neurons and a linear output layer.
We use the $\tanh$ activation function
$\phi_j(\xx)=\tanh(\ba_j\cdot\xx+b_j)$, $j=1,\dots,M$. The parameters $\ba_j\in\R^d$ and
$b_j\in\R$ are sampled independently from $U[-1,1]$ and remain fixed. We use random seed $1$. For each fixed $M$, all methods use the same random
features.

We select the interior collocation points from Cartesian grids and remove a narrow strip near the boundary. Boundary points are placed along polygon edges, on polyhedral faces, or by radial parameterization. The standard RNN method fits the interior and boundary equations with unit weight, whereas the BP method fits only the interior equations. We solve the linear systems in double precision using \texttt{numpy.linalg.lstsq} with \texttt{rcond=None}. Let $N_{\mathrm{test}}$ denote the number of test points. The point counts are listed in Table~\ref{tab:summary}.

We measure accuracy by the discrete relative $L^2$ error
\begin{equation*}
e=\frac{\sqrt{\sum_{i=1}^{N_{\mathrm{test}}}\bigl|\hat u(\xx_i)-u(\xx_i)\bigr|^2}}
{\sqrt{\sum_{i=1}^{N_{\mathrm{test}}}\bigl|u(\xx_i)\bigr|^2}},
\end{equation*}
using the same test points for all methods in each comparison. For polygonal domains, we take $10^4$ uniformly distributed interior test points. For two-dimensional curved domains, we take $10^4$ angular-radial samples: $\theta$ is uniform and $\rho=r(\theta)\sqrt{\xi}$, where $\xi$ is uniform on $[0,1]$. These samples are generally not uniform in area, and the error is computed in the discrete norm on these points. In three dimensions, we use denser filtered grids. The number of test points depends on the domain, as listed in the table.

We refer to the methods based on Wachspress coordinates, mean value coordinates, and radial lifting as the BP-W, BP-M, and BP-R methods, respectively. In three dimensions, these names refer to the corresponding Wachspress, mean value, and spherical radial liftings.
The ratio in the last column of Table~\ref{tab:summary} is $e_{\mathrm{RNN}}/e_{\mathrm{BP}}$, which measures the error reduction of the BP method relative to the standard RNN method.

\begin{table}[tb]
\caption{Relative $L^2$ errors on different domains at $M=400$.}\label{tab:summary}
\centering\scriptsize
\setlength{\tabcolsep}{2.4pt}
\begin{tabular}{lrrrlrrr}
\toprule
Domain & $N_{\mathrm{int}}$ & $N_{\mathrm{bnd}}$ & $N_{\mathrm{test}}$ & Method & RNN & BP & Ratio\\
\midrule
\multicolumn{8}{l}{\textbf{Two-dimensional domains}} \\
Irregular 12-gon & 2432 & 180 & 10000 & BP-W & $1.96\times10^{-9}$ & $4.37\times10^{-12}$ & 448.0 \\
 &  &  &  & BP-M & $1.96\times10^{-9}$ & $6.62\times10^{-12}$ & 296.0 \\
200-gon & 2676 & 185 & 10000 & BP-W & $2.86\times10^{-9}$ & $6.02\times10^{-12}$ & 474.9 \\
 &  &  &  & BP-M & $2.86\times10^{-9}$ & $7.15\times10^{-12}$ & 399.9 \\
Koch & 1892 & 420 & 10000 & BP-M & $1.07\times10^{-9}$ & $2.78\times10^{-12}$ & 383.6 \\
cross & 1248 & 236 & 10000 & BP-M & $9.94\times10^{-10}$ & $2.45\times10^{-12}$ & 405.2 \\
flower & 1536 & 197 & 10000 & BP-R & $9.74\times10^{-9}$ & $8.59\times10^{-10}$ & 11.3 \\
kite & 1994 & 212 & 10000 & BP-R & $8.12\times10^{-8}$ & $6.41\times10^{-9}$ & 12.7 \\
peanut & 2872 & 150 & 10000 & BP-R & $7.77\times10^{-10}$ & $4.15\times10^{-11}$ & 18.7 \\
\midrule
\multicolumn{8}{l}{\textbf{Three-dimensional domains}} \\
Hexagonal prism & 840 & 624 & 8008 & BP-W & $5.94\times10^{-6}$ & $1.75\times10^{-7}$ & 34.0 \\
 &  &  &  & BP-M & $5.94\times10^{-6}$ & $1.01\times10^{-6}$ & 5.9 \\
Dodecahedron & 433 & 304 & 8344 & BP-W & $1.82\times10^{-6}$ & $4.06\times10^{-8}$ & 44.7 \\
 &  &  &  & BP-M & $1.82\times10^{-6}$ & $1.99\times10^{-7}$ & 9.1 \\
Slanted notch & 1502 & 682 & 5115 & BP-M & $3.57\times10^{-5}$ & $2.76\times10^{-6}$ & 13.0 \\
Through hole & 1296 & 864 & 4680 & BP-M & $6.98\times10^{-5}$ & $1.20\times10^{-5}$ & 5.8 \\
Double cavity & 1674 & 816 & 5704 & BP-M & $5.48\times10^{-5}$ & $8.56\times10^{-6}$ & 6.4 \\
peanut & 3296 & 1152 & 7608 & BP-R & $2.35\times10^{-8}$ & $3.94\times10^{-10}$ & 59.5 \\
flower & 2300 & 1152 & 5228 & BP-R & $4.09\times10^{-7}$ & $9.23\times10^{-9}$ & 44.3 \\
egg & 3044 & 1152 & 6924 & BP-R & $5.74\times10^{-8}$ & $7.37\times10^{-10}$ & 77.9 \\
\bottomrule
\end{tabular}
\end{table}

\subsection{Two-dimensional tests}
\label{sec:2d}

In this subsection, we test the BP method on convex and nonconvex polygons and star-shaped curved domains. We compare the errors for different feature counts and numbers of interior collocation points, and report the boundary errors on polygonal domains.

\subsubsection{Polygonal domains}
\label{sec:2dpoly}

\begin{figure}[!htbp]
\centering
\includegraphics[width=0.9\textwidth]{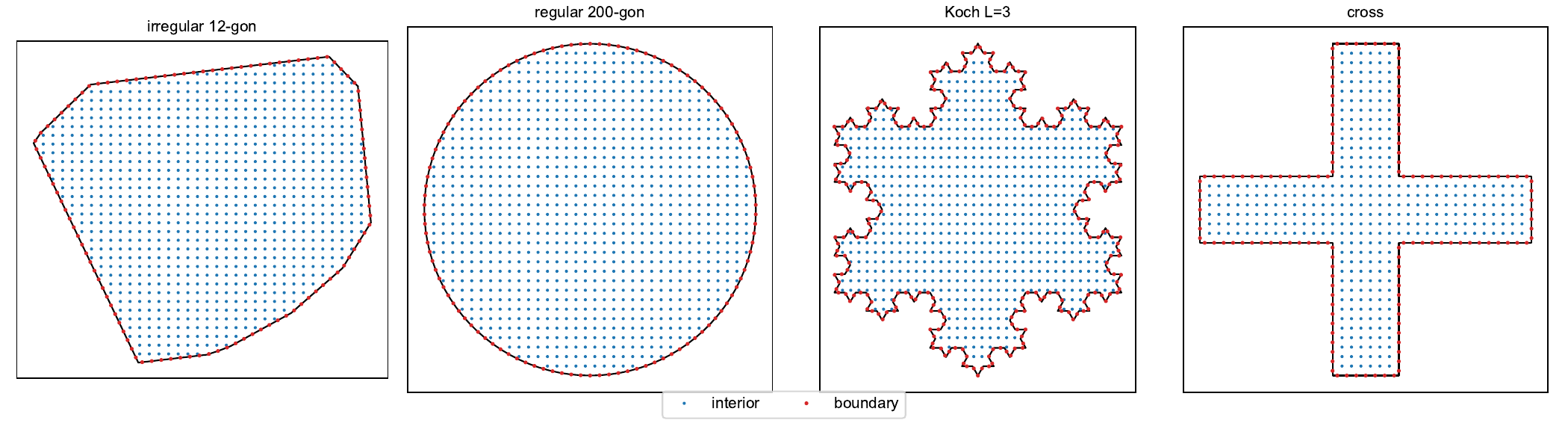}
\caption{Collocation points for the four two-dimensional polygonal test domains.
Blue dots are interior points and red dots are boundary points.}
\label{fig:2dpolycolloc}
\end{figure}

We consider two convex polygons, an irregular polygon with $12$ sides and a regular polygon with $200$ sides,
and two nonconvex polygons, the level-$3$ Koch snowflake and a cross-shaped domain.
We apply both the BP-W and BP-M methods to the convex polygons. For the nonconvex polygons,
we use only the BP-M method because Wachspress coordinates require convexity.
Figure~\ref{fig:2dpolycolloc} shows the collocation points.

\begin{figure}[!htbp]
\centering
\begin{minipage}{0.70\textwidth}
\centering
\begin{subfigure}[t]{0.50\linewidth}
\centering
\includegraphics[width=\linewidth,trim=0bp 0bp 378bp 0bp,clip]{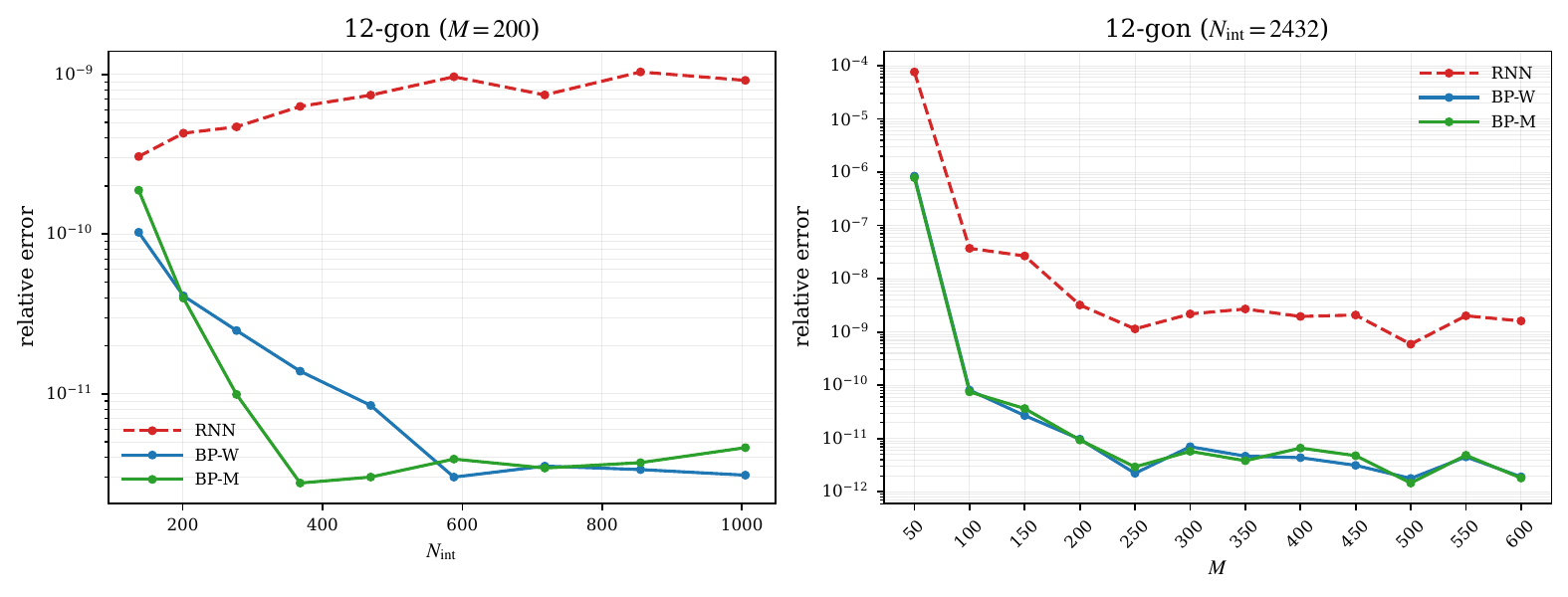}
\caption{Fixed $M=200$; varying $N_{\mathrm{int}}$.}
\label{fig:12gonconv-a}
\end{subfigure}\hfill%
\begin{subfigure}[t]{0.50\linewidth}
\centering
\includegraphics[width=\linewidth,trim=378bp 0bp 0bp 0bp,clip]{figure_02_12gon_conv.pdf}
\caption{Fixed $N_{\mathrm{int}}=2432$; varying $M$.}
\label{fig:12gonconv-b}
\end{subfigure}
\end{minipage}
\caption{Relative $L^2$ errors on the irregular polygon with $12$ sides.}
\label{fig:12gonconv}
\end{figure}

\begin{figure}[!htbp]
\centering
\includegraphics[width=0.35\textwidth]{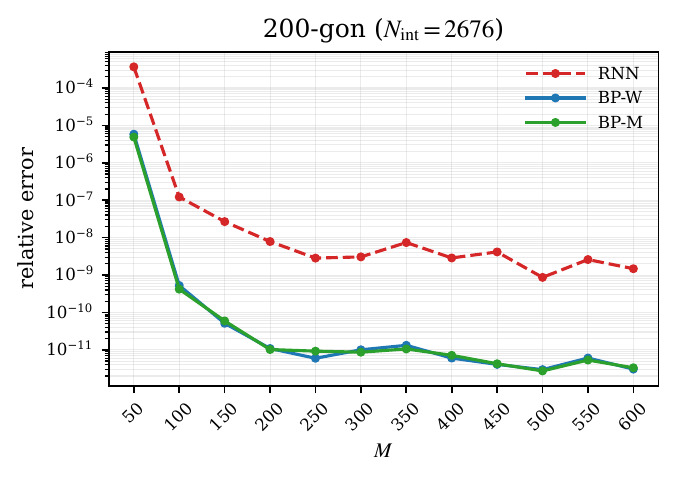}
\caption{Relative $L^2$ errors versus $M$ (fixed $N_{\mathrm{int}}$) on the regular polygon with $200$ sides.}
\label{fig:polyMconv}
\end{figure}

\paragraph{Convex polygons}
Figure~\ref{fig:12gonconv} shows the errors as the number of interior collocation points increases at fixed $M$, and as the number of features increases at fixed sampling. Figure~\ref{fig:polyMconv} shows the errors as $M$ varies on the regular polygon with $200$ sides. For the values of $M$ considered, the BP-W and BP-M methods reduce the relative $L^2$ errors on both convex polygons by factors of approximately $60$--$1000$ compared with the standard RNN method. At fixed $M$, the BP method also gives lower errors as the number of interior collocation points varies.

Table~\ref{tab:2dbnd} gives the maximum absolute error over $10^4$ boundary points distributed by arclength. We use $N_{\mathrm{int}}=855$ for the irregular $12$-sided polygon and $N_{\mathrm{int}}=688$ for the Koch snowflake. The BP method's boundary errors are about $10^{-13}$--$10^{-11}$. Large output coefficients can amplify small floating-point errors, so the measured boundary errors need not be at machine precision.

\begin{figure}[!htbp]
\centering
\includegraphics[width=0.70\textwidth]{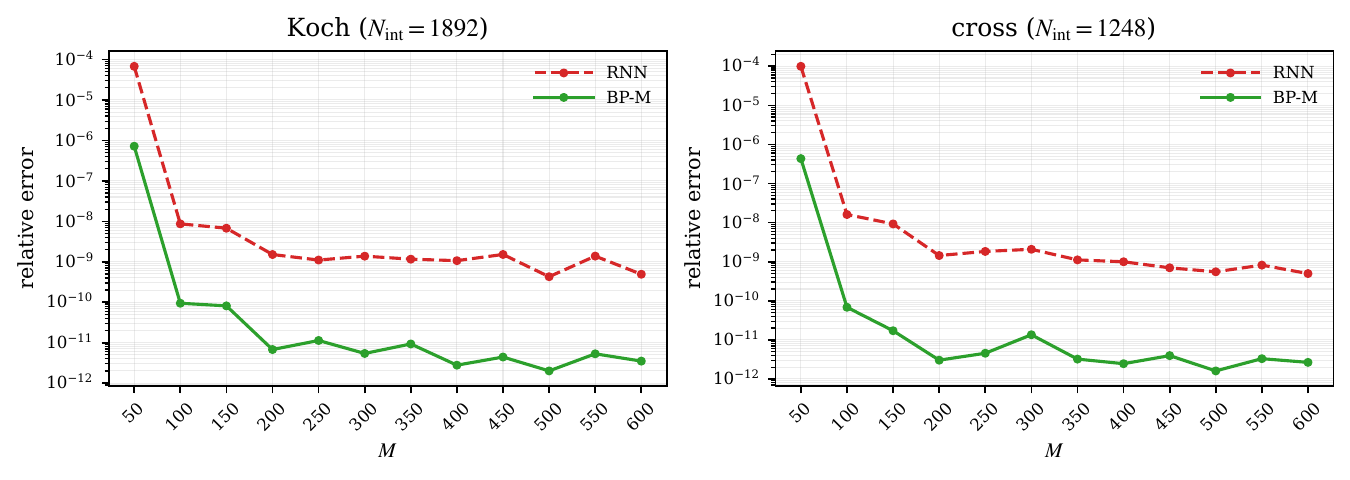}
\caption{Relative $L^2$ errors versus $M$ (fixed $N_{\mathrm{int}}$) on the Koch snowflake and the cross-shaped domain.}
\label{fig:polyMnon}
\end{figure}

\begin{figure}[!htbp]
\centering
\includegraphics[width=\textwidth]{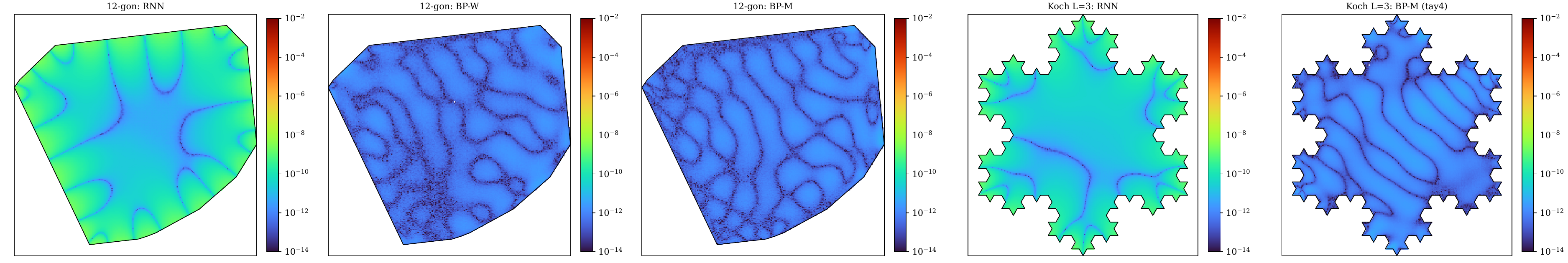}
\caption{Absolute error images of the standard RNN and BP methods on two-dimensional
polygonal regions ($M=400$).}
\label{fig:2derror}
\end{figure}

\paragraph{Nonconvex polygons}
We use the BP-M method for the Koch snowflake and the cross-shaped domain. Since the signed mean value coordinates place some evaluation points beyond the finite edges, we use the quartic boundary-data extension in Section~\ref{sec:extension}. For the feature counts in \rev{Figure}~\ref{fig:polyMnon}, the BP-M method reduces the errors on both domains by factors of approximately \mbox{$80$--$540$} compared with the standard RNN method. The pointwise errors are shown in \rev{Figure}~\ref{fig:2derror}.

\begin{table}[!htbp]
\centering
\caption{Boundary errors of the standard RNN and BP methods on
two-dimensional polygonal domains.}
\label{tab:2dbnd}
\small
\setlength{\tabcolsep}{2.8pt}
\begin{tabular}{l l c c c c}
\toprule
region & method & $M=50$ & $M=100$ & $M=200$ & $M=400$ \\
\midrule
\multirow{3}{*}{12-sided polygon} & RNN & $2.95\times10^{-4}$ & $1.16\times10^{-7}$ & $5.82\times10^{-9}$ & $3.04\times10^{-9}$ \\
 & BP-W & $1.01\times10^{-11}$ & $1.23\times10^{-11}$ & $1.18\times10^{-12}$ & $5.74\times10^{-13}$ \\
 & BP-M & $7.58\times10^{-12}$ & $1.43\times10^{-11}$ & $1.44\times10^{-12}$ & $6.03\times10^{-13}$ \\
\midrule
\multirow{2}{*}{Koch snowflake} & RNN & $3.70\times10^{-4}$ & $7.32\times10^{-8}$ & $1.68\times10^{-9}$ & $2.43\times10^{-9}$ \\
 & BP-M & $4.23\times10^{-12}$ & $1.23\times10^{-11}$ & $7.92\times10^{-13}$ & $2.89\times10^{-13}$ \\
\bottomrule
\end{tabular}
\end{table}

\FloatBarrier
\subsubsection{Star-shaped curved domains}
\label{sec:2dcurved}

\begin{figure}[!htbp]
\centering
\includegraphics[width=0.75\textwidth]{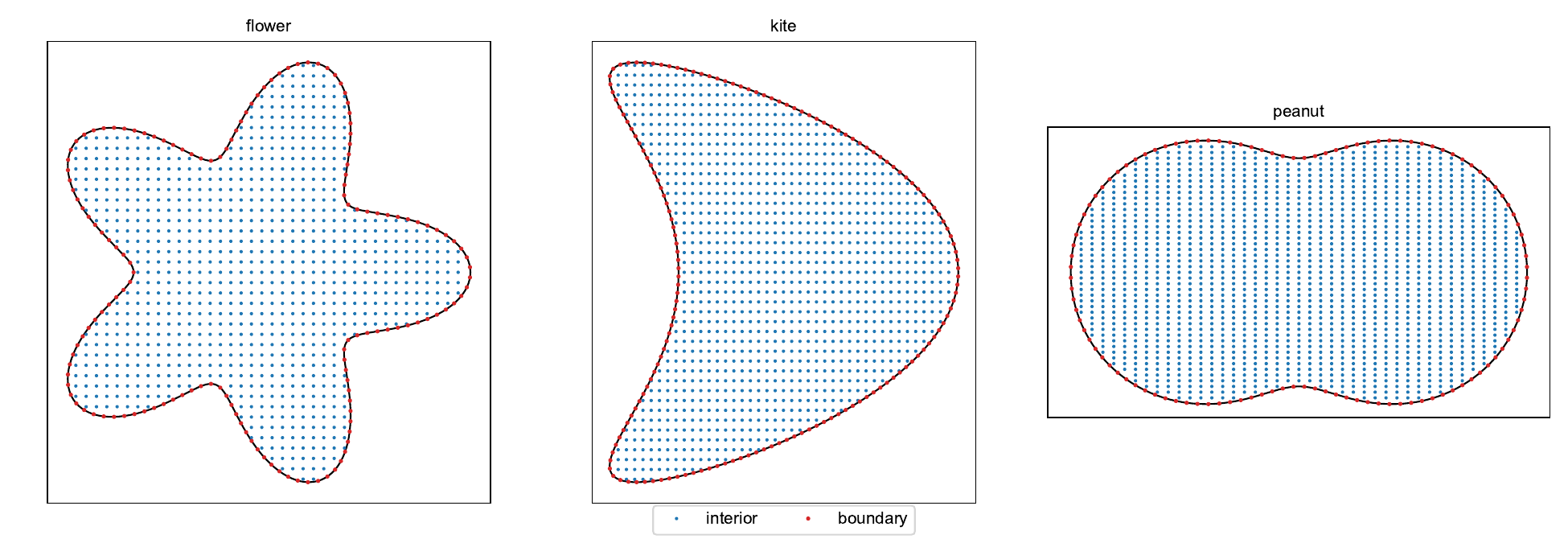}
\caption{Collocation points for the three two-dimensional star-shaped curved test
domains.}
\label{fig:2dcurvedcolloc}
\end{figure}

We apply the BP-R method to three star-shaped curved domains: the five-petal flower-shaped domain, the kite-shaped domain,
and the peanut-shaped domain. Each domain is parameterized radially about an interior center $\mathbf c=(0,0)$.
Their boundary parameterizations are
\begingroup\color{black}
\begin{align*}
\boldsymbol\gamma_{\mathrm{flower}}(t)&=(1+0.3\cos(5t))(\cos t,\sin t),\\
\boldsymbol\gamma_{\mathrm{kite}}(t)&=(\cos t+0.65\cos(2t)-0.65,\,1.5\sin t),\\
\boldsymbol\gamma_{\mathrm{peanut}}(t)&=\sqrt{\cos^2t+\tfrac14\sin^2t}\,(\cos t,\sin t),
\qquad 0\le t<2\pi.
\end{align*}
\endgroup
Figure~\ref{fig:2dcurvedcolloc} shows their collocation points.

\begin{figure}[!htbp]
\centering
\begin{minipage}{0.70\textwidth}
\centering
\begin{subfigure}[t]{0.50\linewidth}
\centering
\includegraphics[width=\linewidth,trim=0bp 0bp 378bp 0bp,clip]{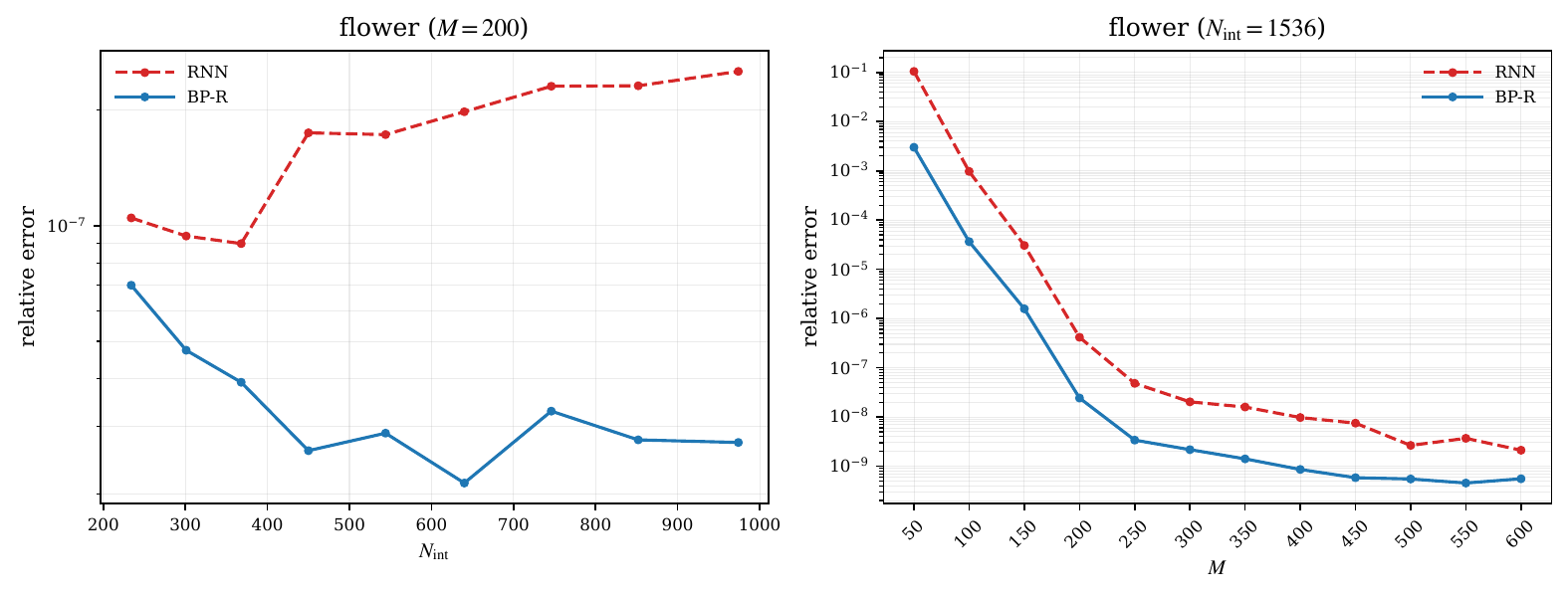}
\caption{Fixed $M=200$; varying $N_{\mathrm{int}}$.}
\label{fig:flowerconv-a}
\end{subfigure}\hfill%
\begin{subfigure}[t]{0.50\linewidth}
\centering
\includegraphics[width=\linewidth,trim=378bp 0bp 0bp 0bp,clip]{figure_07_flower_conv.pdf}
\caption{Fixed $N_{\mathrm{int}}=1536$; varying $M$.}
\label{fig:flowerconv-b}
\end{subfigure}
\end{minipage}
\caption{Relative $L^2$ errors on the flower-shaped domain.}
\label{fig:flowerconv}
\end{figure}

\begin{figure}[!htbp]
\centering
\includegraphics[width=0.70\textwidth]{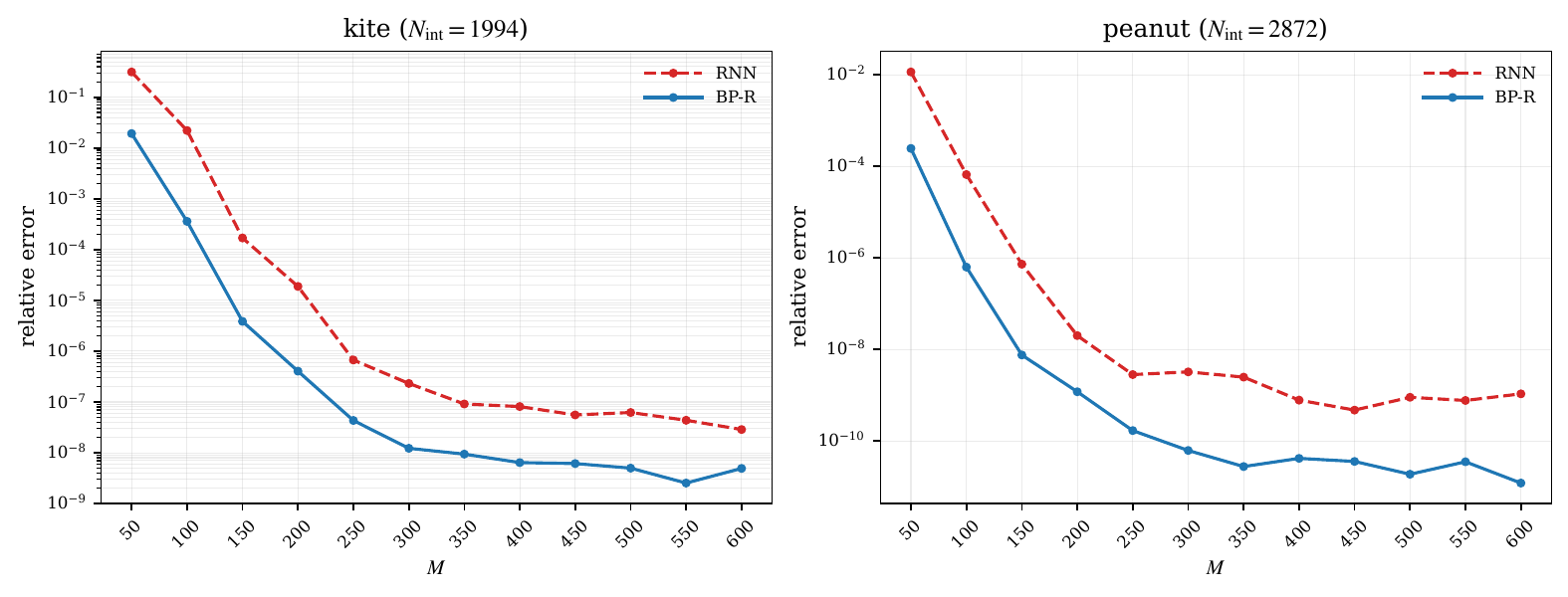}
\caption{Relative $L^2$ errors versus $M$ (fixed $N_{\mathrm{int}}$) on the kite-shaped and peanut-shaped domains.}
\label{fig:kitepeanutM}
\end{figure}

\begin{figure}[!htbp]
\centering
\includegraphics[width=0.75\textwidth]{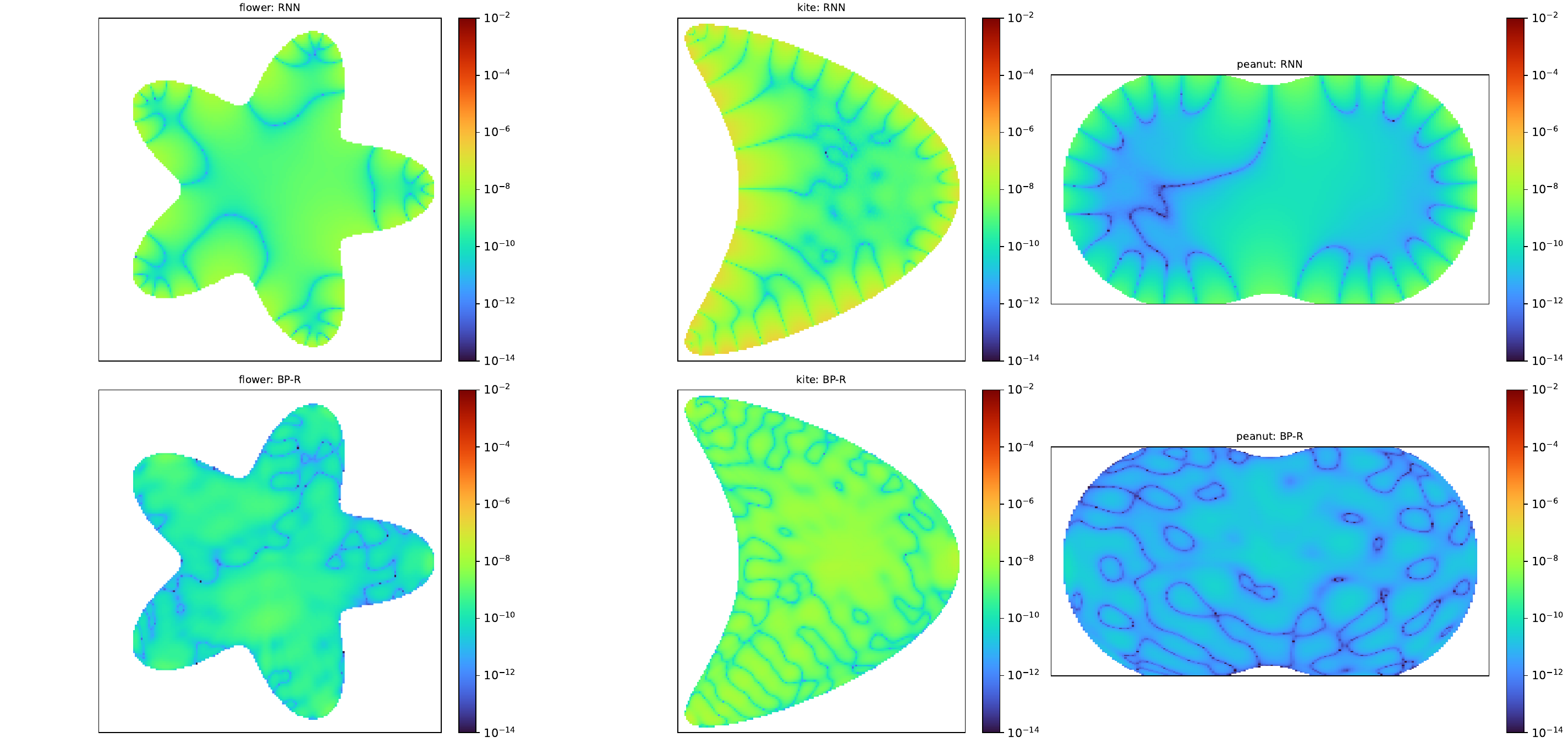}
\caption{Absolute error images of the standard RNN and BP methods on two-dimensional
star-shaped curved domains ($M=400$).}
\label{fig:2dcurvederror}
\end{figure}

Figures~\ref{fig:flowerconv} and \ref{fig:kitepeanutM} show how the errors vary with the numbers of interior collocation points and features. For the values of $M$ considered, the BP-R method reduces the errors on all three domains by factors of approximately $4$--$105$ compared with the standard RNN method. For the flower-shaped domain, the BP-R method also gives lower errors for different numbers of interior collocation points at fixed $M$. The pointwise errors are presented in \rev{Figure}~\ref{fig:2dcurvederror}.

\FloatBarrier
\subsection{Three-dimensional tests}
\label{sec:3d}

In this subsection, we test the BP methods on convex and nonconvex polyhedra and star-shaped curved domains. We compare their errors with those of the standard RNN method for different feature counts.

\subsubsection{Polyhedral domains}
\label{sec:3dpoly}

\begin{figure}[!htbp]
\centering
\includegraphics[width=0.23\textwidth]{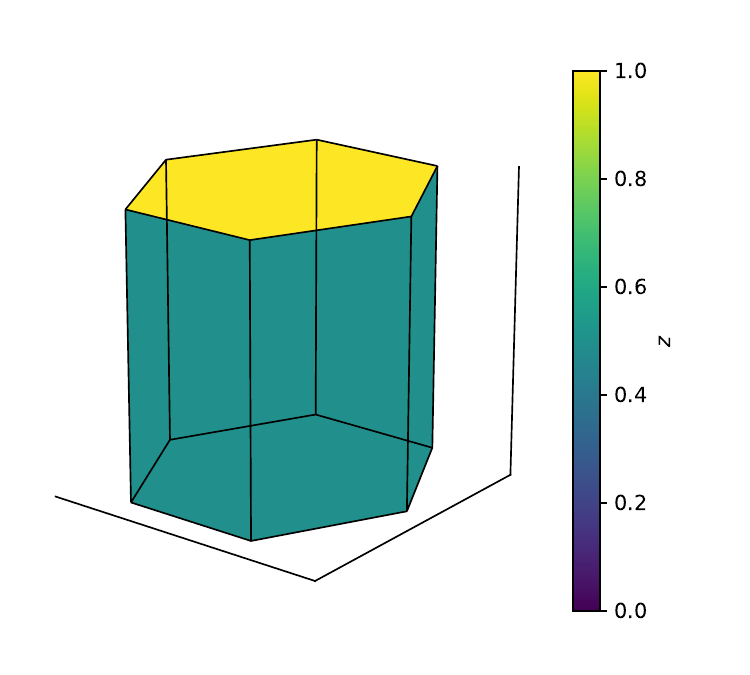}\hfill
\includegraphics[width=0.23\textwidth]{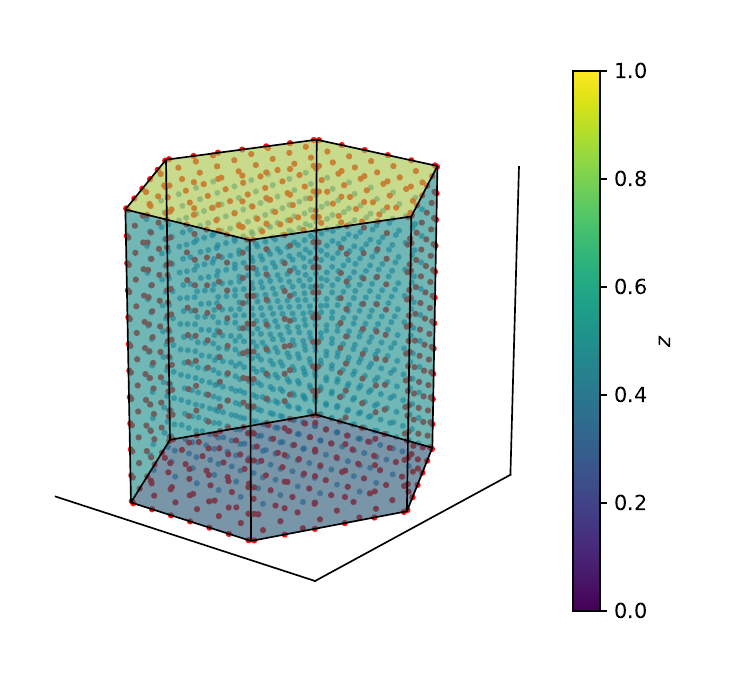}\hfill
\includegraphics[width=0.23\textwidth]{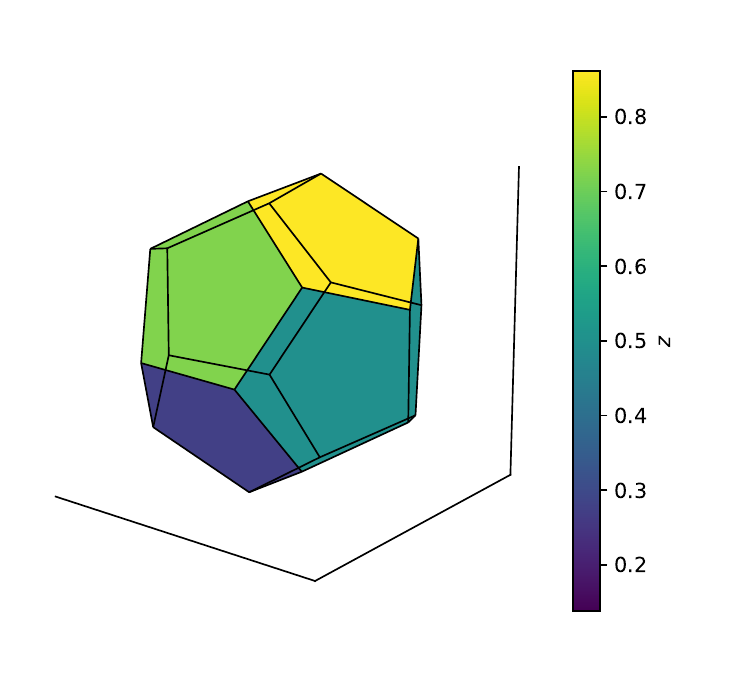}\hfill
\includegraphics[width=0.23\textwidth]{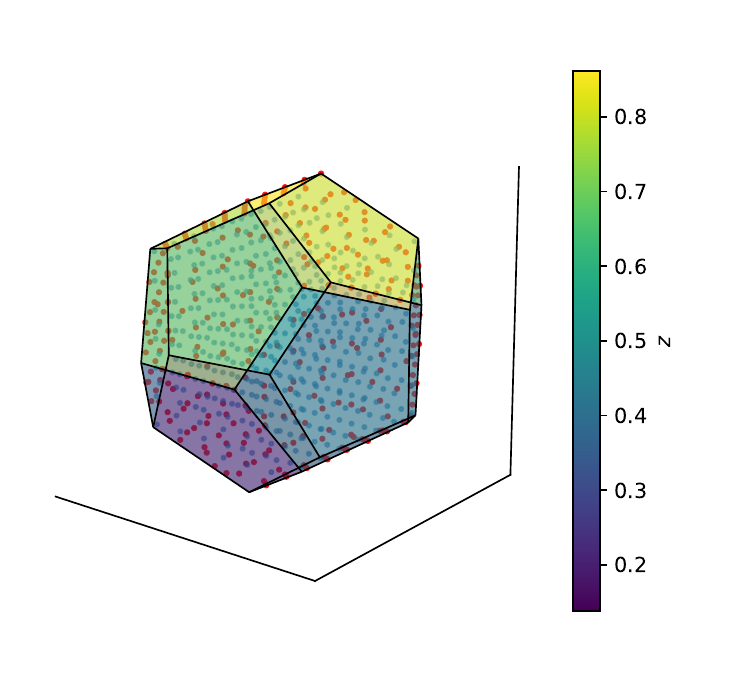}
\caption{From left to right: the hexagonal prism region and its collocation points,
the regular dodecahedron region and its collocation points.}
\label{fig:3dpolycolloc}
\end{figure}

\begin{figure}[!htbp]\centering
\includegraphics[width=0.70\textwidth]{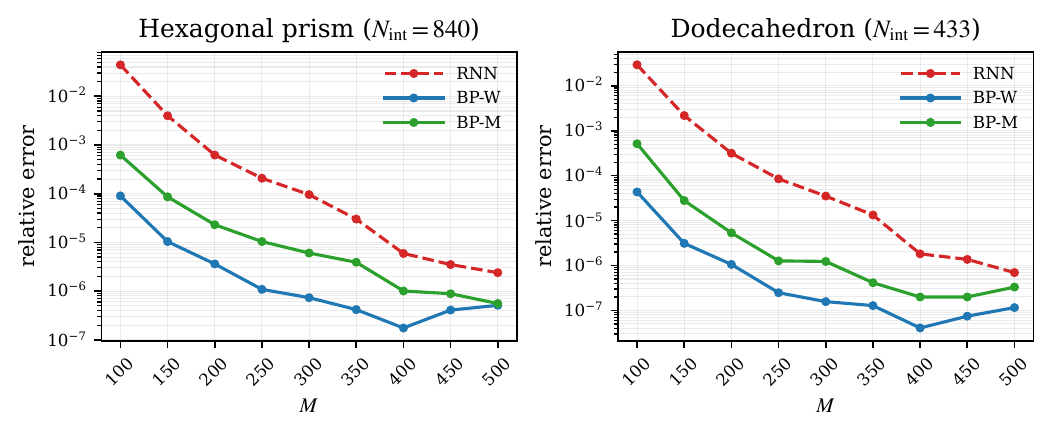}
\caption{Three-dimensional convex polyhedra: relative $L^2$ errors versus $M$ at fixed interior sampling. Left: hexagonal prism; right: dodecahedron.}\label{fig:3dhex}\label{fig:3ddodeca}
\end{figure}

\paragraph{Convex polyhedra}
We first compare the standard RNN method with the BP-W and BP-M methods on a regular hexagonal prism and a regular
dodecahedron. The prism has unit height and a regular hexagonal base with
circumradius $0.5$. The dodecahedron is scaled so that its bounding box is
$[0,1]^3$. Figure~\ref{fig:3dpolycolloc} shows the domains and their collocation
points. The errors for different feature counts are given in \rev{Figure}~\ref{fig:3dhex}. Over the range shown,
the BP-W and BP-M methods reduce the errors on both domains by factors of approximately
$5$--$700$ and $2$--$80$, respectively, compared with the standard RNN method. The BP-W method gives the larger
error reductions.

\begin{figure}[!htbp]
\centering
\includegraphics[width=0.25\textwidth]{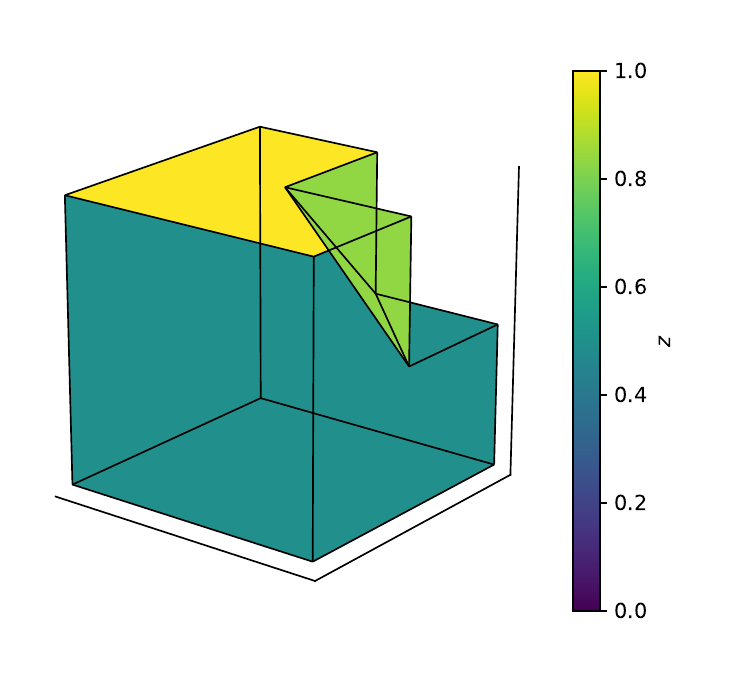}\hspace{0.025\textwidth}%
\includegraphics[width=0.25\textwidth]{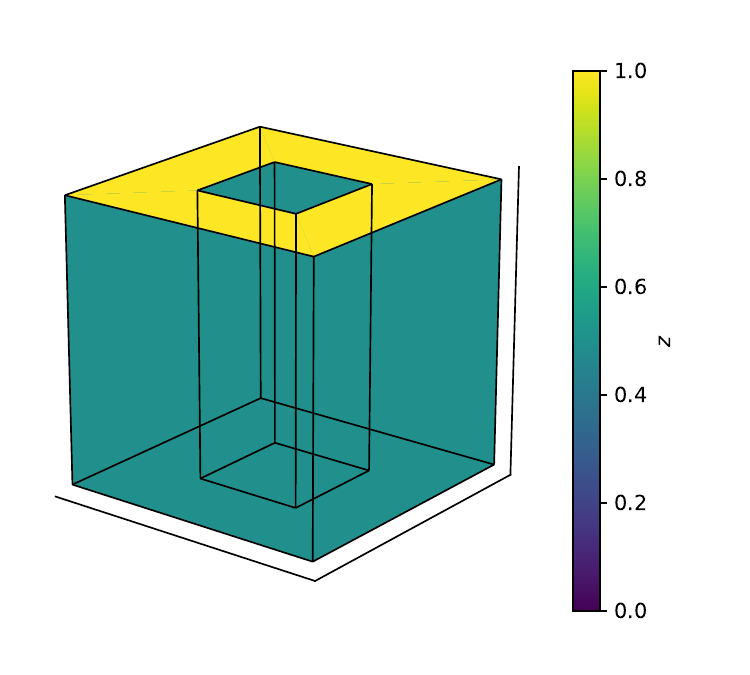}\hspace{0.025\textwidth}%
\includegraphics[width=0.25\textwidth]{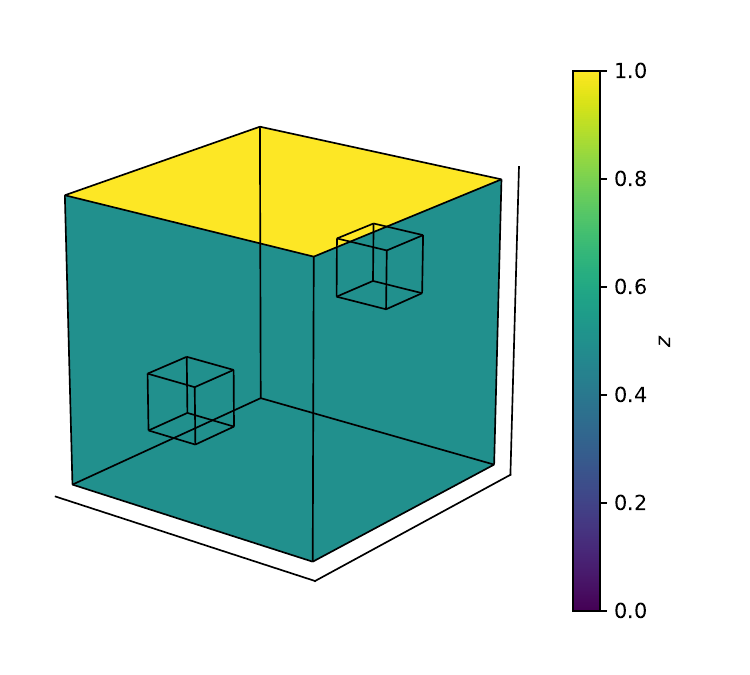}
\caption{The three nonconvex polyhedra: slanted notch, through hole, and double
cavity.}
\label{fig:3dnonconvex}
\end{figure}

\begin{figure}[!htbp]
\centering
\includegraphics[width=0.90\textwidth]{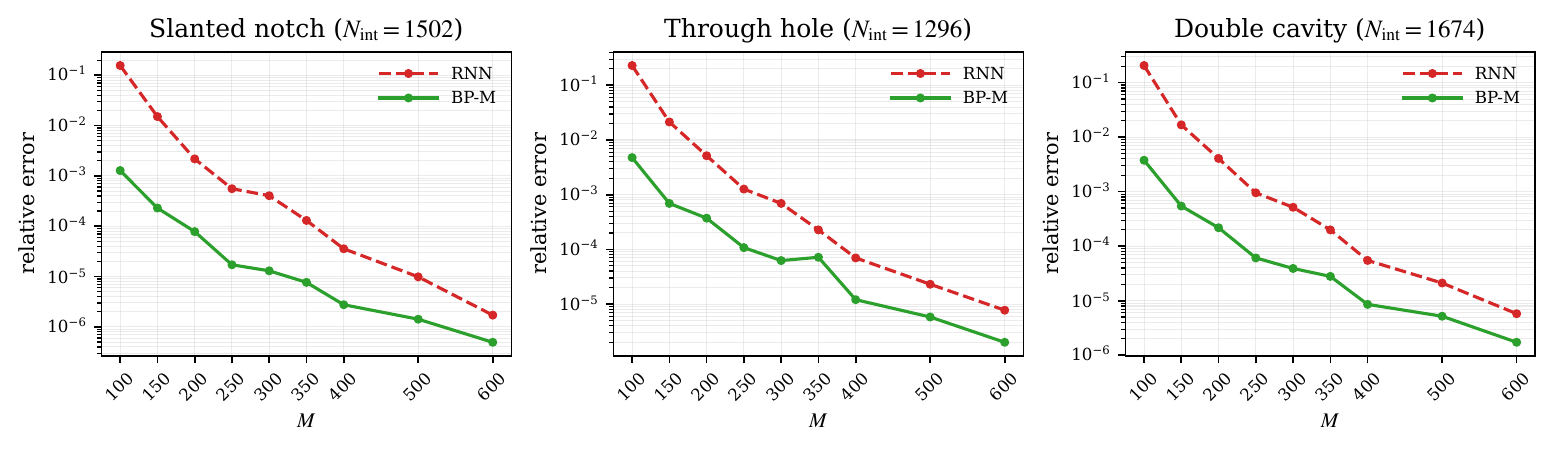}
\caption{Relative $L^2$ errors versus $M$ (fixed $N_{\mathrm{int}}$) on the three
nonconvex polyhedra.}
\label{fig:3dnonconvexM}
\end{figure}

\paragraph{Nonconvex polyhedra}
We next consider the slanted-notch, through-hole, and double-cavity blocks.
The slanted-notch block is obtained by removing the reentrant corner
$\{x,y,z\ge0.5,\ x+y+z\ge2.0\}$ from the unit cube. The through-hole block
removes $[0.3,0.7]^2\times[0,1]$, and the double-cavity block removes
$[0.18,0.38]^3$ and $[0.62,0.82]^3$. The domains and their collocation points are
shown in \rev{Figure}~\ref{fig:3dnonconvex}. The relative $L^2$ errors as $M$ varies at fixed
$N_{\mathrm{int}}$ are reported in \rev{Figure}~\ref{fig:3dnonconvexM}.
We use the quartic face/edge boundary-data extension in Section~\ref{sec:extension}.
For the feature counts shown, the BP-M method reduces the errors on all three domains by
factors of approximately $3$--$120$ compared with the standard RNN method.

\begin{figure}[!htbp]
\centering
\includegraphics[width=0.25\textwidth]{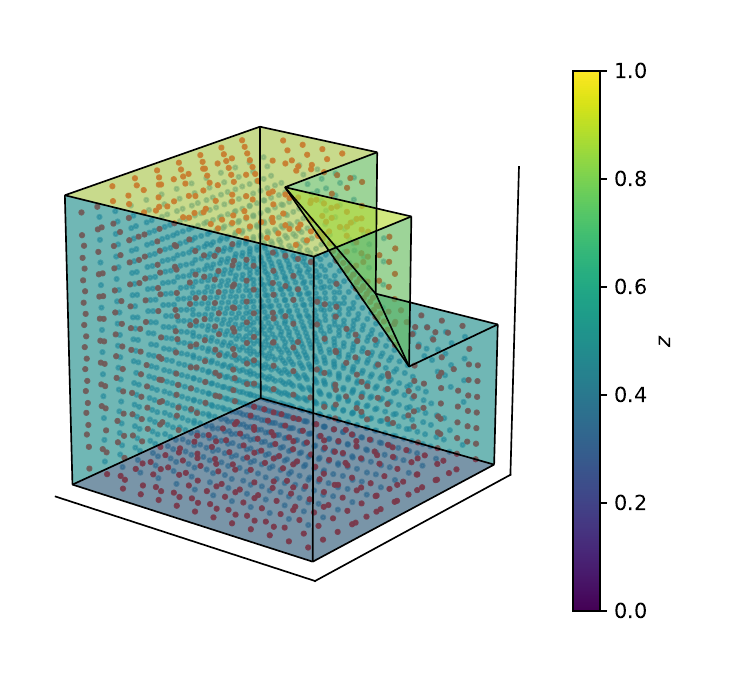}\hspace{0.025\textwidth}%
\includegraphics[width=0.25\textwidth]{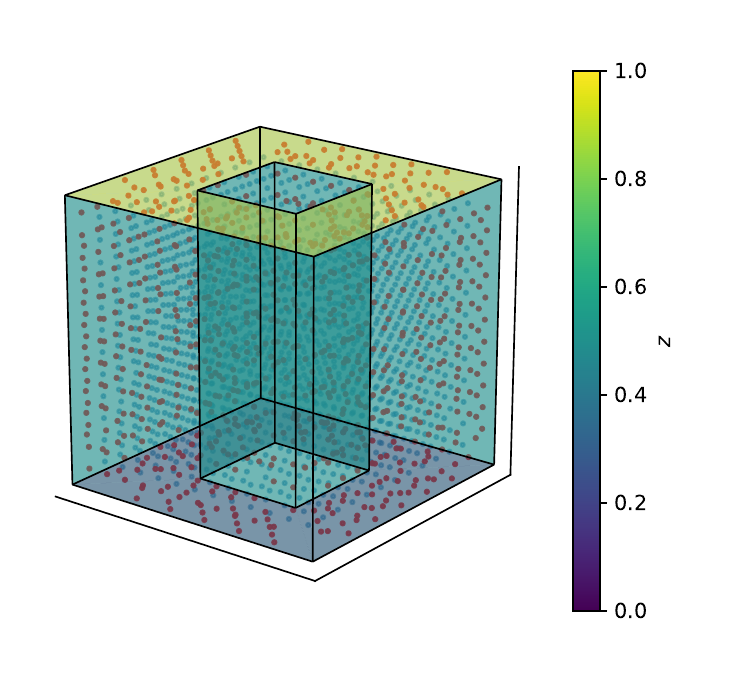}\hspace{0.025\textwidth}%
\includegraphics[width=0.25\textwidth]{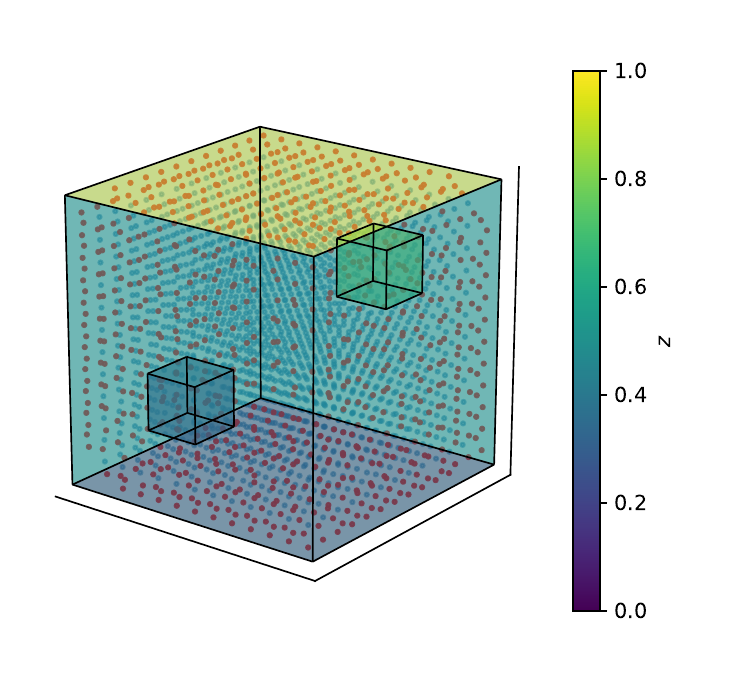}
\caption{Collocation points for the three nonconvex polyhedra.}
\label{fig:3dnonconvexcolloc}
\end{figure}

\FloatBarrier
\subsubsection{Curved domains}
\label{sec:3dcurved}

\begin{figure}[!htbp]
\centering
\includegraphics[width=0.25\textwidth]{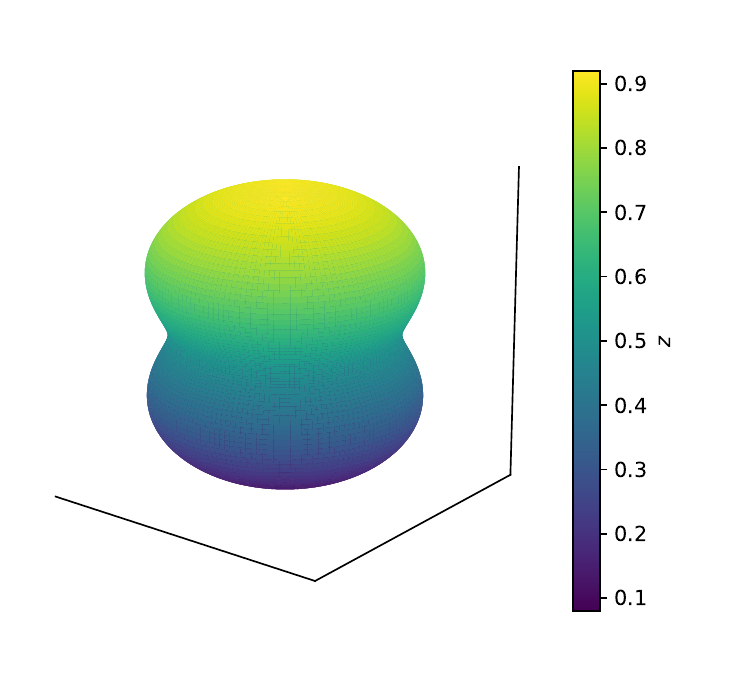}\hspace{0.025\textwidth}%
\includegraphics[width=0.25\textwidth]{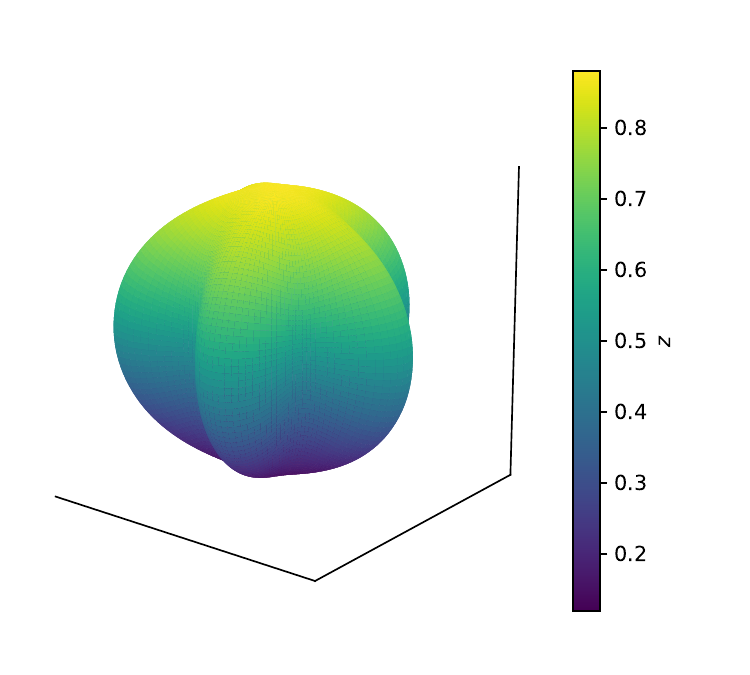}\hspace{0.025\textwidth}%
\includegraphics[width=0.25\textwidth]{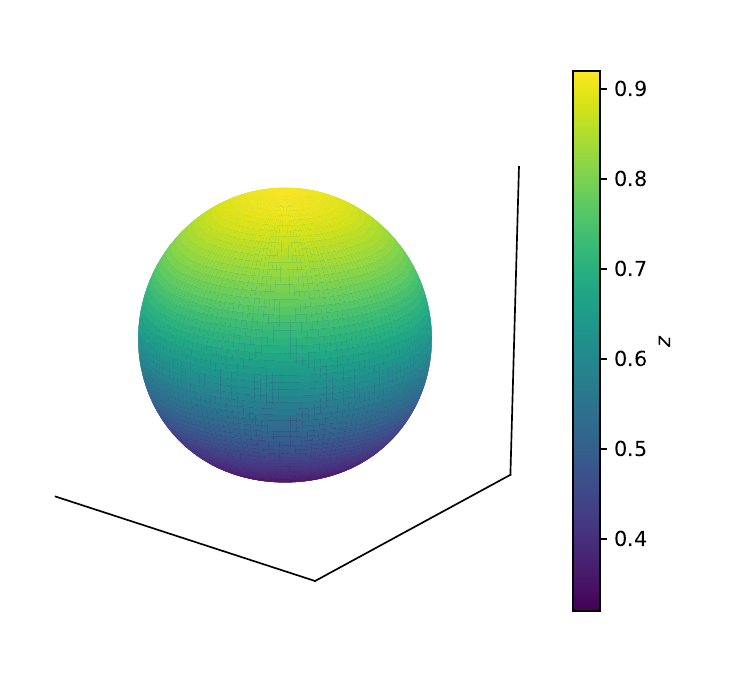}
\caption{The peanut-shaped, flower-shaped, and egg-shaped domains.}
\label{fig:3dcurvedregion}
\end{figure}

\begin{figure}[!htbp]
\centering
\includegraphics[width=0.25\textwidth]{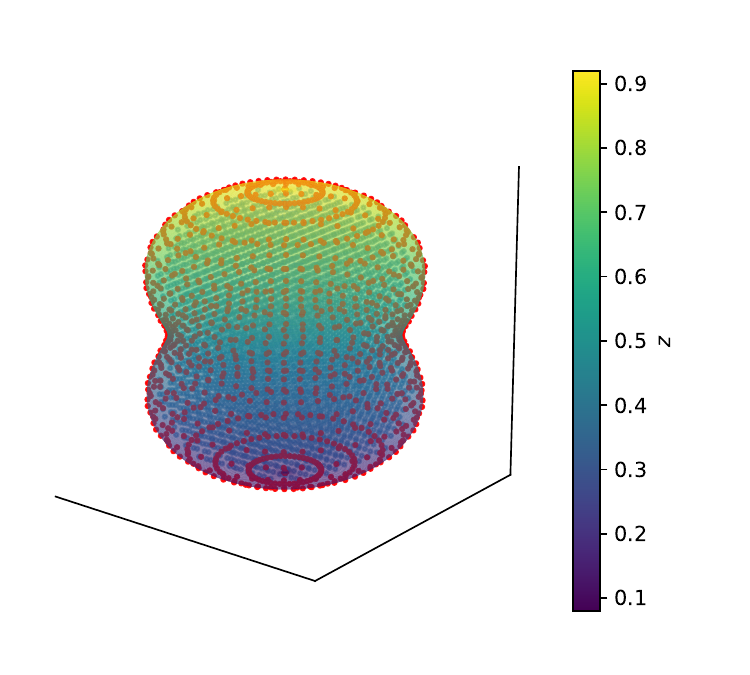}\hspace{0.025\textwidth}%
\includegraphics[width=0.25\textwidth]{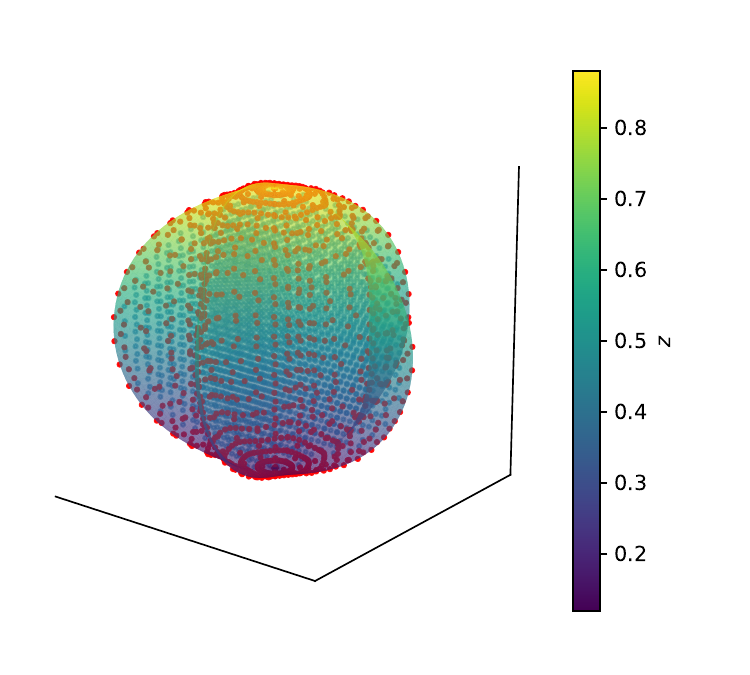}\hspace{0.025\textwidth}%
\includegraphics[width=0.25\textwidth]{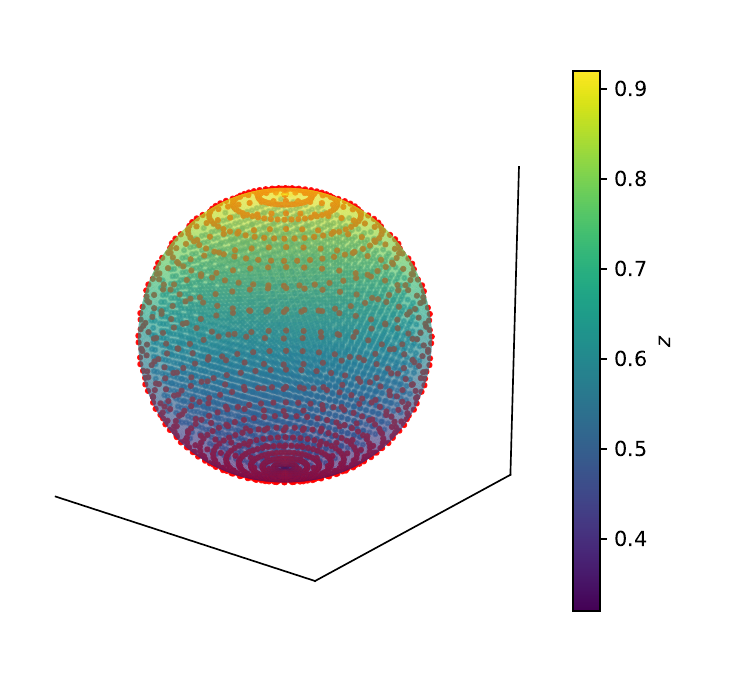}
\caption{Collocation points for the three curved domains.}
\label{fig:3dcurvedcolloc}
\end{figure}

\begin{figure}[!htbp]
\centering
\includegraphics[width=0.30\textwidth]{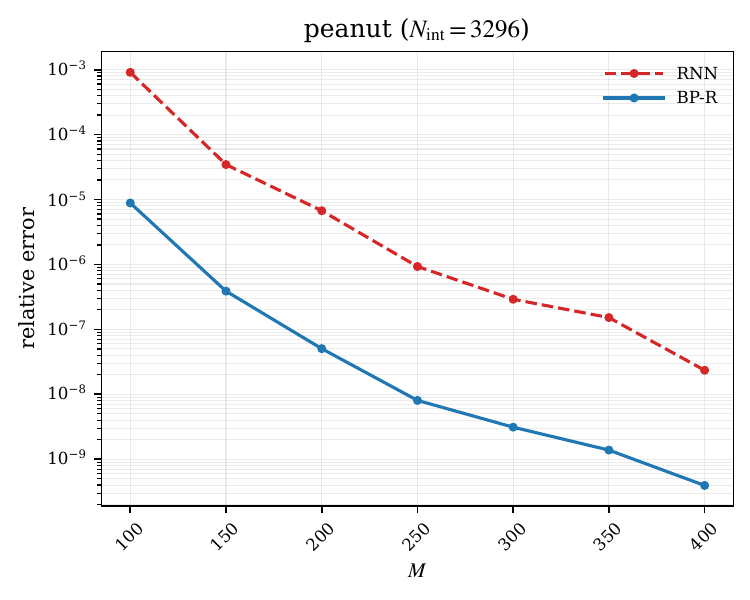}\hfill
\includegraphics[width=0.30\textwidth]{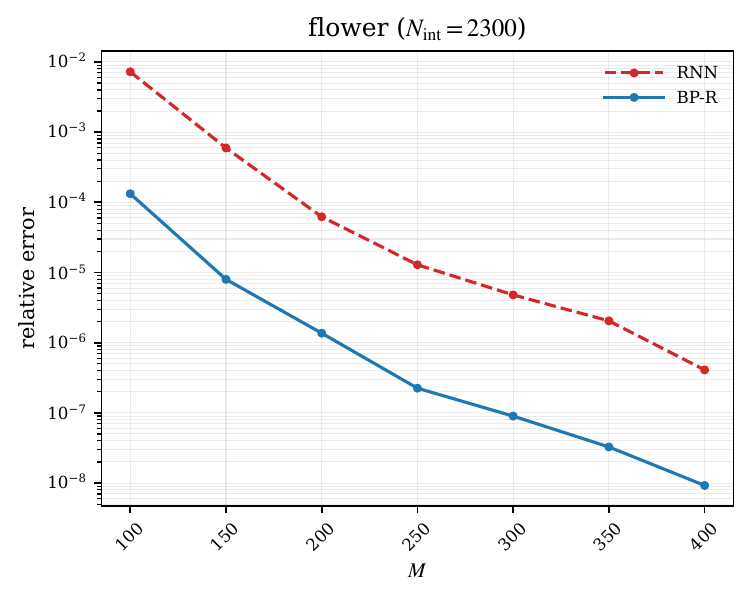}\hfill
\includegraphics[width=0.30\textwidth]{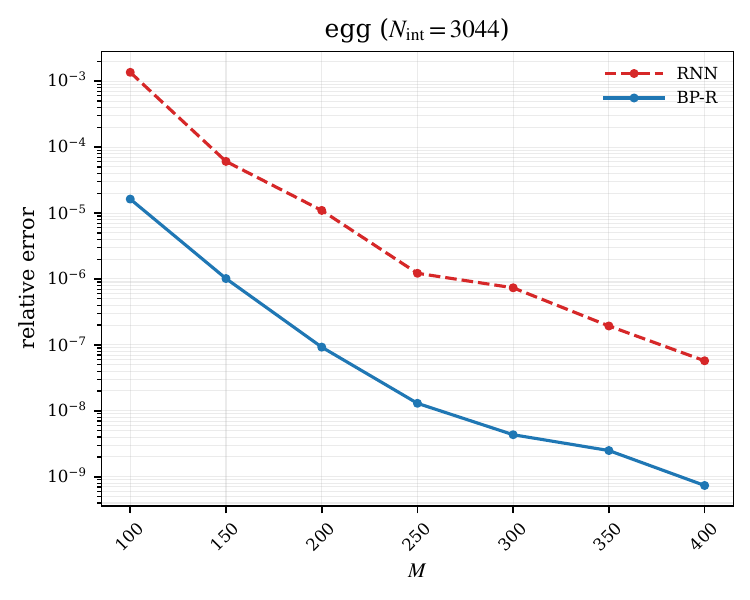}
\caption{Relative $L^2$ errors versus $M$ (fixed $N_{\mathrm{int}}$) on the peanut-shaped,
flower-shaped, and egg-shaped domains.}
\label{fig:3dcurvedM}
\end{figure}

We use the BP-R method on three star-shaped curved domains with center $(0.5,0.5,0.5)$: the peanut-shaped domain
$r=0.30+0.12\cos(2\theta)$ (nonconvex), the five-petal flower-shaped domain
$r=0.38+0.10\sin^3\theta\cos(5\phi)$ (nonconvex), and the egg-shaped domain
$r=0.30+0.12\cos\theta$ (convex). Figures~\ref{fig:3dcurvedregion} and
\ref{fig:3dcurvedcolloc} show the domains and their collocation points, respectively.
The relative $L^2$ errors as $M$ varies at fixed $N_{\mathrm{int}}$ are given in
\rev{Figure}~\ref{fig:3dcurvedM}. For the feature counts shown, the BP-R method gives lower relative
$L^2$ errors on all three domains, with reductions of approximately $40$--$170$
times compared with the standard RNN method.

\FloatBarrier
\section{Conclusions}
\label{sec:conclusion}

In this paper, we propose a boundary-preserving randomized neural network method for two- and three-dimensional \rev{Poisson equations with Dirichlet boundary conditions}. The Dirichlet boundary conditions are satisfied by construction, so no penalty weights need to be balanced between the interior and boundary residuals. The method treats different geometries within the same boundary-preserving framework. We construct geometry-dependent liftings for convex and nonconvex polygons and polyhedra, and for star-shaped curved domains satisfying the required radial representation. In the numerical experiments, the BP method reduces the relative $L^2$ errors by factors of approximately $1.5$--$1000$ compared with the standard RNN method. For most two-dimensional \rev{polygon cases}, the errors are reduced by more than two orders of magnitude. The three-dimensional curved-domain tests also give reductions ranging from \rev{$40$--$170$} times.

\FloatBarrier
\section*{CRediT authorship contribution statement}
Huifang Zhou: Conceptualization, Methodology, Software, Validation, Formal analysis, Investigation, Visualization, Writing -- original draft, Writing -- review and editing.

\section*{Declaration of competing interest}
The author declares no known competing financial interests or personal relationships that could have appeared to influence the work reported in this paper.

\section*{Funding}
This work was supported by the Science Challenge Project (Grant No. TZ2026004) and the National Natural Science Foundation of China (Grant No. 12671473).

\section*{Data availability}
The code and numerical data supporting the findings of this study are available from the corresponding author upon reasonable request.

\end{document}